\documentclass[10pt]{amsart}
\usepackage{bbm}
\usepackage{mathrsfs}
\usepackage{amsfonts} 
\usepackage{graphicx,latexsym,bm,amsmath,amssymb,verbatim,multicol,lscape}
\newtheorem{thm}{Theorem}[section]

\theoremstyle{definition}

\theoremstyle{remark}

\theoremstyle{problem}
\newtheorem{pro}[thm]{Problem}
\numberwithin{equation}{section}

\begin{document}
\title[The LCM of consecutive quadratic progression terms]
{The least common multiple of consecutive quadratic progression terms}%
\author{Shaofang Hong}
\address{Mathematical College, Sichuan University, Chengdu 610064, P.R. China}
\email{sfhong@scu.edu.cn, s-f.hong@tom.com, hongsf02@yahoo.com}
\author{Guoyou Qian}
\address{Center for Combinatorics, Nankai University, Tianjin 300071, P.R. China}
\email{qiangy1230@163.com, qiangy1230@gmail.com}
\thanks{The work was supported partially by National Science Foundation of
China Grant \#11371260, by the Ph.D. Programs
Foundation of Ministry of Education of China Grant \#20100181110073 and
by the Postdoctoral Science Foundation of China Grant \#2013M530109}
\date{\today}%
\begin{abstract}
Let $k$ be an arbitrary given positive integer and let $f(x)\in
{\mathbb Z}[x]$ be a quadratic polynomial with $a$ and $D$ as its
leading coefficient and discriminant, respectively. Associated to
the least common multiple ${\rm lcm}_{0\le i\le k}\{f(n+i)\}$
of any $k+1$ consecutive terms in the quadratic progression $\{f(n)\}_{n\in
\mathbb{N}^*}$, we define the function $g_{k, f}(n):=(\prod_{i=0}^{k}|f(n+i)|)
/{\rm lcm}_{0\le i\le k}\{f(n+i)\}$ for all integers $n\in \mathbb{N}^*\setminus
Z_{k, f}$, where $Z_{k,f}:=\bigcup_{i=0}^k\{n\in \mathbb{N}^*: f(n+i)=0\}$.
In this paper, we first show that $g_{k,f}$ is eventually periodic
if and only if $D\ne a^2i^2$ for all integers $i$ with $1\le i\le k$.
Consequently, we develop a detailed $p$-adic analysis of $g_{k, f}$ and
determine its smallest period. Finally, we obtain asymptotic formulas of
$\log {\rm lcm}_{0\le i\le k}\{f(n+i)\}$ for all quadratic polynomials
$f$ as $n$ goes to infinity.
\end{abstract}

\maketitle

\section{\bf Introduction and the statements of the main results}

The study of least common multiple of consecutive positive
integers was initiated by Chebyshev \cite{[Ch]} for the first
significant attempt to prove prime number theorem. Motivated by
Chebyshev's work, one naturally expects to investigate the least
common multiple of consecutive terms in any given sequence of
positive integers. For the least common multiple of the first $n$
terms of a given sequence of positive integers, some results were
obtained by several authors. Hanson \cite{[Ha]} and Nair
\cite{[N]} got the upper and lower bound of ${\rm lcm}_{1\le
i\le n}\{i\}$ respectively. Bateman, Kalb and Stenger \cite{[BKS]}
obtained an asymptotic formula for the least common multiple of
arithmetic progressions. Farhi \cite{[F1]} \cite{[F2]} and Farhi
and Kane \cite{[FK]} studied the least common multiple of some
finite sequences of integers. Hong and Feng \cite{[HF]}, Hong and
Yang \cite{[HY2]}, Kane and Kominers \cite{[KK]} and Wu, Tan and
Hong \cite{[WTH]} investigated the least common multiple of finite
arithmetic progressions. Hong, Qian and Tan \cite{[HQT]} obtained
an asymptotic estimate for the least common multiple of a
sequence of products of linear polynomials. Qian and Hong
\cite{[QH]} got an asymptotic formula for the least common
multiple of consecutive arithmetic progression terms.

The investigation of periodic arithmetic functions has been a common topic in
number theory for a long time. The readers are refereed to \cite{[A]} and
\cite{[M]} for the related background information.
When investigating the least common multiple ${\rm lcm}_{0\le i\le k}$
$\{n+i\}$ of any $k+1$ consecutive integers with $k$ being a fixed
positive integer, Farhi \cite{[F2]} introduced the arithmetic function
$\bar g_k$ defined for positive integer $n$ by
$\bar g_k(n):=\frac{\prod_{i=0}^k(n+i)}{{\rm lcm}_{0\le i\le k}\{n+i\}}.$
Farhi showed that $\bar g_k$ is periodic with $k!$ as its period. Let $\bar P_k$
be the smallest period of $\bar g_k$. Then  $\bar P_k|k!$. At the end
of \cite{[F2]}, Farhi  posed the open problem of determining the
smallest period $\bar P_k$. Hong and Yang \cite{[HY]} improved the
period $k!$ to ${\rm lcm}_{1\le i\le k}\{i\}$ and proposed a
conjecture stating  that $\frac{{\rm lcm}_{1\le i\le
k+1}\{i\}}{k+1}$ divides $\bar P_k$. Farhi and Kane \cite{[FK]}
proved the Hong-Yang conjecture and finally determined the exact
value of $\bar P_k$. Throughout, let ${\mathbb Q}$, ${\mathbb Z}$ and
${\mathbb N}$ denote the field of rational numbers, the ring of integers
and the set of nonnegative integers, respectively.
Define ${\mathbb N}^*:={\mathbb N}\setminus\{0\}$. Let $b\in
\mathbb{N}$ and $a, k\in {\mathbb N}^*$. Define $L_k:={\rm lcm}_{1\le
i\le k}\{i\}.$ Hong and Qian \cite{[HQ]} defined the arithmetic function
$g_{k,a,b}:\mathbb{N}^*\longrightarrow \mathbb{N}^*$ by
$
g_{k,a,b}(n):=\frac{\prod_{i=0}^{k}(b+(n+i)a)} {{\rm lcm}_{0\le i\le
k}\{b+(n+i)a\}},
$
and proved that $g_{k,a,b}$ is periodic and determined
the exact value of the smallest period of $g_{k,a,b}$.
Let $f(x)$ be a quadratic polynomial with integer coefficients.
Associated to ${\rm lcm}_{0\le i\le k}\{f(n+i)\}$, we define
the function $g_{k, f}$ for all positive integers
$n\in \mathbb{N}^*\setminus Z_{k, f}$ by
\begin{align}\label{1.1}
g_{k, f}(n):=\frac{\prod_{i=0}^k |f(n+i)|}{{\rm lcm}_{0\le i\le
k}\{f(n+i)\}},
\end{align}
where $Z_{k,f}:=\bigcup_{i=0}^k\{n\in \mathbb{N}^*: f(n+i)=0\}.$
Recall that an arithmetic function $\alpha $ is called
{\it eventually periodic} if there are positive integers $t_0$
and $n_0$ such that $\alpha (n+t_0)=\alpha (n)$ holds for all
integers $n\ge n_0$. If $f(x)=x^2+1$, then it is proved in
\cite{[QTH]} that $g_{k, f}$ is periodic and also the smallest period
is determined there. We can easily find that there are
quadratic polynomials $f$ such that $g_{k, f}$ is not eventually periodic.
Therefore one naturally asks the following interesting question.

\begin{pro}\label{pro1.1} \cite {[Ho]} Let $f(x)$ be a polynomial of
degree two with integer coefficients and $g_{k,f}$ be defined as in (\ref{1.1}).
Characterize $f(x)$ such that $g_{k,f}$ is eventually periodic. If $g_{k,f}$ is
eventually periodic, what is the smallest period of $g_{k,f}$?
\end{pro}
In this paper, we mainly focus on the least common multiple of any $k+1$ consecutive
terms in the quadratic progression. Our main goal is to study Problem \ref{pro1.1}.
First we use a well-known identity of Hua \cite{[Hu]} to determine all the quadratic
integer polynomials $f$ such that $g_{k, f}$ is eventually periodic. Then we transfer
the computation of the smallest period of $g_{k, f}$ into the computation of
the local smallest periods. To determine the local smallest periods, we will first investigate
the minimal distance among the roots of quadratic congruences and then develop
a detailed local analysis. To state the main results, let's introduce and recall
some notation. As usual, for any prime number $p$, we let $\nu _{p}$ be the normalized
$p$-adic valuation of ${\mathbb Q}$, i.e., $\nu _p(s)=t$ if
$p^{t}\parallel s$. Let ${\rm gcd}(s, t)$ denote the greatest
common divisor of any integers $s$ and $t$. For any real number $x$,
by $\lfloor x\rfloor$ we denote the largest integer no more than
$x$. Throughout this paper, we always let $a\ge 1$ and $f(x)=ax^2+bx+c$ be any given
quadratic polynomial with integer coefficients, and let
$D:=b^2-4ac$ be the discriminant of $f$. Define
$D_4:=\frac{D}{4^{\lfloor \frac{\nu _2(D)}{2}\rfloor}}$ and $D_p:=\frac{D}{p^{\nu _p(D)}}$
for any odd prime $p$. Then $D_4$ is equal to $\frac{D}{2^{\nu _2(D)}}$
if $\nu _2(D)$ is even, and equals $\frac{2D}{2^{\nu _2(D)}}$
if $\nu _2(D)$ is odd. So $D_4\not\equiv 0\pmod 4$,
$\nu _2(D)$ and $D_4$ hold the reverse parity. Let
$\big(\frac{\cdot }{p}\big)$ denote the Legendre symbol. For any
positive integer $k$, we define $B_k:={\rm lcm}_{1\le i\le k}\{i(a^2i^2-D)\}$
and
\begin{align}\label{1.2}
A_k:=\frac{B_k}{\xi_2\Big(\displaystyle\prod_{p\ne 2,
p|\gcd(a, b)}p^{\nu _p(B_k)}\Big)\Big(\prod_{p\nmid 2aD,
(\frac{D}{p})=-1}p^{\nu _p(B_k)}\Big)\Big(\prod_{p\nmid 2a, p|D}\eta_p\Big)},
\end{align}
where

\begin{align}\label{1.3}
\xi_2={\left\{
  \begin{array}{rl}
  1, & \mbox{if}\ 2|a, 2\nmid b\ \mbox{and}\ \nu _2(k+1)<\nu _2(B_k),\\
  2^{2\nu _2(L_k)}, & \mbox{if}\ 2\nmid a,\  k<2^{\lfloor \frac{\nu _2(D)}
  {2}\rfloor} \ \mbox{and}\ \nu _2(k+1)<\nu _2(L_k),\\
  2^{\nu _2(B_k)-\lfloor \frac{\nu _2(D)}
  {2}\rfloor}, & \mbox{if} \  2\nmid a, k\ge
  2^{\lfloor \frac{\nu _2(D)}{2}\rfloor}, D_4\not\equiv1\pmod8
  \ \mbox{and}\ \nu _2(k+1)<\lfloor \frac{\nu _2(D)}{2}\rfloor,\\
2^{\nu _2(D)+1},&  \mbox{if} \ 2\nmid a, k\ge
  2^{\lfloor \frac{\nu _2(D)}{2}\rfloor} \ \mbox{and}\
  D_4\equiv1 \pmod 8,\\
2^{\nu _2(B_k)}, & \mbox{otherwise}\\
 \end{array}
\right.}
\end{align}
and
\begin{align}\label{1.4}
\eta_p={\left\{
  \begin{array}{rl}
  p^{2\nu _p(L_k)}, & \mbox{if} \ k<p^{\lceil \frac{\nu _p(D)}
  {2}\rceil} \ \mbox{and}\ \nu _p(k+1)<\nu _p(L_k),\\
  p^{\nu _p(B_k)-\lceil \frac{\nu _p(D)}
  {2}\rceil}, & \mbox{if} \
  k\ge p^{\lceil \frac{\nu _p(D)}{2}\rceil},\
  \nu _p(k+1)<\lceil \frac{\nu _p(D)}{2}\rceil \\ &
  \mbox{and either} \ 2\nmid \nu _p(D) \  \mbox{or}\ (\frac{D_p}{p})=-1,\\
p^{\nu _p(D)},&   \mbox{if} \  k\ge
  p^{\lceil \frac{\nu _p(D)}{2}\rceil}, \nu _p(k+1)<\nu _p(B_k)-\nu _p(D),\\
  & 2|\nu _p(D) \ \mbox{and}\ (\frac{D_p}{p})=1,\\
p^{\nu _p(B_k)}, & \mbox{otherwise.}\\
 \end{array}
\right.}
\end{align}
We can now state the first main result of this paper as follows.

\begin{thm}\label{thm1.1} Let $f(x)$ be a polynomial of degree two
with integer coefficients and $k$ be a positive integer. Let
$g_{k,f}$ be defined as in (\ref{1.1}). Then
$g_{k,f}$ is eventually periodic if and only if
$D\ne a^2i^2$ for all integers $i$ with $1\le i\le k$.
If $g_{k,f}$ is eventually periodic, then its smallest
period is equal to $A_k$ except that $\nu _p(k+1)\ge \nu _p(A_k)\ge 1$
for at most one odd prime $p$ such that either $p|a$ and $p\nmid b$
or $p\nmid 2aD$ and $(\frac{D}{p})=1$, in which case its smallest
period equals $A_k/p^{\nu _p(A_k)}$.
\end{thm}
\noindent Therefore Theorem \ref{thm1.1} answers completely Problem \ref{pro1.1}.

By \cite{[HQ]}, we know that for any integers $a\ge 1$ and $b$,
one has $\log {\rm lcm}_{0\le i\le k}\{a(n+i)+b\}\sim (k+1) \log
n$ as $n$ goes to infinity. Now using Theorem \ref{thm1.1}, we can deduce
the following asymptotic estimates of $\log{\rm lcm}_{0\le i\le
k}\{f(n+i)\}$ for all quadratic polynomials $f$ with
integer coefficients as $n$ tends to infinity. This is the second
main result of this paper.
\begin{thm}\label{thm1.2}
Let $k$ be a positive integer. Then
$\log {\rm lcm}_{0\le i\le k}\{f(n+i)\}=C_{k, D}\log n+o(\log n)$,
where $C_{k, D}=2(k+1)$ if $D\ne a^2i^2$ for all integers $i$ with
$1\le i\le k$, and $C_{k, D}=k+i_0+1$ if $D=a^2i_0^2$ for some
integer $i_0$ with $1\le i_0\le k$.
\end{thm}

This paper is organized as follows. We first show in Section 2
that $g_{k,f}$ is eventually periodic if and only if
$D\ne a^2i^2$ for all integers $i$ with $1\le i\le k$.
Subsequently, we give a formula which factors the global smallest
period $P_{k, f}$ into the product of the local smallest
periods $P_{p, k, f}$. Then in Section 3, we study
the structure of the roots of quadratic congruences and introduce
the concept of the minimal distance among the roots of quadratic
congruences. Consequently, we develop some arithmetic properties
of the minimal distance. In Section 4, we provide a detailed
$p$-adic analysis of $g_{k,f}$, and then using the arithmetic
results obtained in Section 3, we arrive at explicit formulas
of the local snallest periods $P_{p, k, f}$. In Section 5,
by using the results presented in Section 4, we show Theorem
\ref{thm1.1}. Some examples are also given in Section 5 to
demonstrate the validity of Theorem \ref{thm1.1}. Finally,
in Section 6, we give the proof of Theorem \ref{thm1.2}.

\section{\bf A characterization on $f$ such that $g_{k,f}$ is eventually periodic}

Throughout this section, we always let $f$ be a quadratic primitive
polynomial with integer coefficients (i.e. the greatest common divisor
of all the coefficients of $f$ is 1). We first characterize all the
quadratic primitive polynomials $f$ with integer coefficients such that
$g_{k,f}$ is eventually periodic. Subsequently, we transfer
the smallest period problem into a local analysis problem.
We begin with the following result which answers the first
part of Problem \ref{pro1.1}.

\noindent{\bf Theorem 2.1.} {\it Let $k$ be a positive integer.
The function $g_{k,f}$ is eventually periodic if and only if
$D\ne a^2i^2$ for all integers $i$ with $1\le i\le k$.
Furthermore, if $g_{k,f}$ is eventually periodic, then $B_k$ is its period.}
\begin{proof}
First we show the necessity part. Let $g_{k,f}$ be eventually periodic.
Suppose that $D=a^2i_0^2$ for some integer $i_0$
with $1\le i_0\le k$. Then $f(x)$ is reducible. We may let
$f(x)=ax^2+bx+c:=(a_1x+b_1)(a_2x+b_2)$
with $\gcd(a_1, b_1)=\gcd(a_2, b_2)=1$ and $a_1, a_2\in \mathbb{N}^*$. Then
$D=(a_2b_1-a_1b_2)^2=a_1^2a_2^2i_0^2.$
In other words, we have
$a_2b_1-a_1b_2=\pm a_1a_2i_0$ and so $a_2(b_1\pm a_1i_0)=a_1b_2$. It
follows that $a_1=a_2$ and $b_2=b_1\pm a_1i_0$. So we can write $f$ as
$f(x)=(a_1x+b_1)(a_1x+b_1\pm a_1i_0).$

If $b_2=b_1+a_1i_0$, then for any positive integer $n$, $a_1n+b_1+a_1i_0$
divides $\frac{|f(n)f(n+i_0)|}{{\rm lcm}(f(n), f(n+i_0))}.$
Hence $(a_1n+b_1+a_1i_0)\big| g_{k, f}(n)$ which implies that
$g_{k, f}(n)\ge a_1n+b_1+a_1i_0$.

If $b_2=b_1-a_1i_0$, we obtain that $(a_1n+b_1)\big |g_{k, f}(n)$ and
$g_{k, f}(n)\ge a_1n+b_1$. Thus $g_{k, f}(n)$ tends to infinity as
$n$ tends to infinity. This is impossible since $g_{k, f}$ is eventually
periodic implies that $g_{k, f}(n)$ is bounded. The necessity part is proved.

Consequently, we show the sufficiency part. Let $D\ne a^2i^2$
for all integers $i$ with $1\le i\le k$. Then
we have $B_k={\rm lcm}_{1\le i\le k}\{i(a^2i^2-D)\}\ne 0$. For any
given positive integer $n\in \mathbb{N}^*\setminus Z_{k,f}$,
we derive from the identity
$(2an+3aj-ai+b)f(n+i)-(2an+3ai-aj+b)f(n+j)=(j-i)(a^2(j-i)^2-D)$
that
\begin{align}\label{3.1}
\gcd\big( f(n+i), f(n+j)\big)|{\rm lcm }_{0\le i<j\le
k}\{(j-i)(a^2(j-i)^2-D)\}=B_k.
\end{align}
It follows that
$
\gcd\big( f(n+i), f(n+j)\big)|f(n+i\pm B_k)\ \mbox{and}\
\gcd\big(f(n+i), f(n+j)\big)|f(n+j\pm B_k).
$
Thus we have
$
\gcd\big(f(n+i), f(n+j)\big)| \gcd\big(f(n+B_k+i), f(n+B_k+j)\big)
$
and
\begin{align}\label{3.2}
\gcd\big(f(n+i), f(n+j)\big)| \gcd\big(f(n-B_k+i), f(n-B_k+j)\big).
\end{align}
Replacing $n$ by $n+B_k$ in (\ref{3.2}), one gets
$
\gcd\big( f(n+B_k+i), f(n+B_k+j)\big)|\gcd\big( f(n+i), f(n+j)\big).
$
Therefore ${\rm gcd}(f(n+i), f(n+j))={\rm
gcd}(f(n+i+B_k), f(n+j+B_k))$ for any positive integer $n\in
\mathbb{N}^*\setminus Z_{k,f}$ and any integers $i,j$ with $0\le
i<j\le k$. But Theorem 7.3 in Chapter 1 of
\cite{[Hu]} (see Page 11 of \cite{[Hu]}) tells us that
$$g_{k,f}(n)=\prod_{r=1}^k\prod_{0\le i_0<...<i_r\le
k}\big(\gcd\big(f(n+i_0), ..., f(n+i_r) \big)\big)^{(-1)^{r-1}}$$
and
$$
g_{k,f}(n+B_k)=\prod_{r=1}^k\prod_{0\le i_0<...<i_r\le
k}\big(\gcd\big(f(n+B_k+i_0), ..., f(n+B_k+i_r)
\big)\big)^{(-1)^{r-1}}.
$$
Thus $g_{k,f}(n+B_k)=g_{k,f}(n)$ for any positive integer $n\in
\mathbb{N}^*\setminus Z_{k,f}$, and $g_{k,f}$ is eventually
periodic with $B_k$ as its period. The proof of Theorem 2.1 is complete.
\end{proof}

It should be pointed out that the condition that $D\ne a^2i^2$
for all integers $i$ with $1\le i\le k$ is equivalent to saying that
$f(x)$ and $f(x+i)$ are relatively prime in ${\mathbb Q}[x]$ for all
integers $i$ with $1\le i\le k$.
So Theorem 2.1 tells us that $g_{k, f}$ is eventually periodic
if and only if $f(x)$ and ${\rm lcm}_{1\le i\le k}\{f(x+i)\}$
are relatively prime in ${\mathbb Q}[x]$.
After giving a characterization on $f$ so that $g_{k,f}$ is
eventually periodic, we now turn our attention to
determining the smallest period of $g_{k, f}$. Our basic idea is to
transfer the problem into a local analysis problem such that we can
provide a local analysis to $g_{k,f}$. It is clear that only when $f$ is
reducible, $Z_{k,f}$ is not empty. In this case, for all $n\in
Z_{k,f}$, we can always find a positive integer $a_0$ such that
$n+a_0B_k\in \mathbb{N}^*\backslash Z_{k,f}$. By defining
$g_{k,f}(n):=g_{k,f}(n+a_0B_k)$ for each $n\in Z_{k,f}$, we get the
extended periodic arithmetic function
$g_{k,f}:\mathbb{N}^*\longrightarrow \mathbb{N}^*$ with $B_k$ as its
period. In what follows, when mentioning
$g_{k,f}$, we will mean the extended periodic arithmetic function
$g_{k,f}$. For any given prime $p$, we define the arithmetic function
$g_{p,k,f}$ for any positive integer $n$ by
$g_{p,k,f}(n):=v_{p}(g_{k,f}(n))$. If $g_{k,f}$ is periodic with
$P_{k,f}$ as its smallest period, then $g_{p,k,f}$ is periodic and
$P_{k,f}$ is a period of $g_{p,k,f}$. Let $P_{p,k,f}$ be the
smallest period of $g_{p,k,f}$. Then $P_{p, k, f}|P_{k, f}$.
The following result factors the global smallest period $P_{k, f}$
into the product of the local smallest periods $P_{p, k, f}$.

\noindent{\bf Lemma 2.2.} {\it For any prime $p$, $P_{p,k,f}$
divides $p^{\nu _p(B_k)}$. Further, $P_{k,f}=\prod_{p| B_k}P_{p,k,f}.$}
\begin{proof} For any positive integer $n$ and any two
integers $i, j$ with $0\le i<j\le k$, we get by (\ref{3.1}) that
$v_{p}({\rm gcd} (f(n+i), f(n+j)))=\min\{\nu _p(f(n+i)), \nu _p(f(n+j))\}\le
\nu _p(B_k),$
which means that $\nu _p(f(n+i))\le \nu _p(B_k)$ or
$\nu _p(f(n+j))\le \nu _p(B_k)$. Therefore
$\nu _p(f(n+i))\le \nu _p(f(n+i\pm p^{\nu _p(B_k)}))$
or $\nu _p(f(n+j))\le \nu _p(f(n+j\pm p^{\nu _p(B_k)})).$
Hence we derive that
\begin{align*}
 \nu _p\big({\rm gcd} (f(n+i), f(n+j))\big) &\le \min\big(\nu _p
(f(n+i+p^{\nu _p(B_k)})), \nu _p(
f(n+j+p^{\nu _p(B_k)}))\big )\\
 & = v_{p}\big({\rm gcd}
(f(n+i+p^{\nu _p(B_k)}), f(n+j+p^{\nu _p(B_k)}))\big)
\end{align*}
and
\begin{align}\label{3.3}
\nu _p\big({\rm gcd} (f(n+i), f(n+j))\big)\le \nu _p\big({\rm gcd}
(f(n+i-p^{\nu _p(B_k)}), f(n+j-p^{\nu _p(B_k)}))\big).
\end{align}
Replacing $n$ by $n+p^{\nu _p(B_k)}$ in (\ref{3.3}), we obtain
$
\nu _p\big(\gcd (f(n+i+p^{\nu _p(B_k)}), f(n+j+p^{\nu _p(B_k)}))\big)\le
\nu _p\big(\gcd(f(n+i), f(n+j))\big).
$
Thus we have that
$
v_{p}({\rm gcd} (f(n+i+p^{\nu _p(B_k)}), f(n+j+p^{\nu _p(B_k)})))
=v_{p}({\rm gcd} (f(n+i), f(n+j)))
$
for any positive integer $n$ and any two integers $0\le i<j\le k$.
It then follows from Theorem 7.3 in Chapter 1 of \cite{[Hu]} that
$g_{p, k, f}(n)=g_{p,k,f}(n+p^{\nu _p(B_k)})$ for any positive integer
$n$.  Hence we get that $p^{\nu _p(B_k)}$ is a period of $g_{p,k,f}$ and
thus $P_{p,k,f}| p^{\nu _p(B_k)}$. It tells us that $P_{p,k,f}$ are
relatively prime for different prime numbers $p$ and $P_{p,k,f}=1$
for all primes $p\nmid B_k$.

On the other hand, since $P_{q,k,f}|P_{k, f}$ for each
prime $q$, we have
$\prod_{{\rm prime}\ q|B_k}P_{q,k,f}\Big| P_{k, f}.$
Since $P_{q, k,f}=1$ for all primes $q\nmid R_k$, we have
for each prime $p$ and any positive integer $n$ that
$\nu _p(g_{k,f}(n+\prod_{{\rm prime}\ q|B_k}P_{q,k,f}))=\nu _p(g_{k,f}(n)),$
which implies that $\prod_{p|B_k}P_{p,k,f}$ is a period of $g_{k, f}$
and $P_{k, f}\Big|\prod_{p|B_k}P_{p,k,f}.$
Thus the desired result follows immediately. So Lemma 2.2 is proved.
\end{proof}

By Lemma 2.2, to get the global smallest period $P_{k,f}$, it is
sufficient to determine the exact value of the local smallest
period $P_{p,k,f}$ for all the primes $p|B_k$. The remaining part of
the paper will devote to establishing such a local analysis.

\section{\bf Minimal distance among the roots of a quadratic congruence}

Throughout this section, we let $f(x)=ax^2+bx+c$ be an arbitrary primitive
quadratic polynomial with integer coefficients
and let $p$ denote a prime. A natural question is to determine the
roots of the congruence $f(x)\equiv0\pmod {p^e}$ and to investigate
the relation among  distinct roots. Note that the number of roots of
the congruence $x^2\equiv n\pmod{p^e}$ is given in \cite{[Hu]},
where $e$ and $n$ are positive integers such that $p\nmid n$.
Also notice that the problem of distribution of roots of
quadratic congruences to prime modulus  was  investigated by  Duke,
Friedlander, Iwaniec \cite{[DFI]}. Our concern here is the structure
of the roots of the congruence $f(x)\equiv0\pmod {p^e}$.

For any given nonnegative integer $e$, by $S(f, p^e)$ we denote the
set of solutions $x$ with $1\le x\le p^e$ of the congruence
$f(x)\equiv 0\pmod{p^e}$.  Evidently, $S(f, p^0)=\{1\}$.
Throughout, for any $x\in\mathbb{Z}_p$, the ring of
$p$-adic integers, by $\langle x\rangle _{p^e}$ we mean an integer
between 1 and $p^e$ such that $\langle x\rangle _{p^e}\equiv x\pmod {p^e}$.
We begin with the following lemma.\\

\noindent{\bf Lemma 3.1.} {\it Let $e$ be a positive integer and let $p$
be a prime such that $p|a$. Then $S(f, p^e)$ is empty if $p|b$, and equals
$\{\langle s_p\rangle_{p^e}\}$ if $p\nmid b$, where $s_p$ is the unique
solution of the equation $f(x)=0$ in the ring $\mathbb{Z}_p$ of
$p$-adic integers.}
\begin{proof}
If $p|a$ and $p|b$, then $p\nmid c$ since $\gcd(a, b, c)=1$.
Hence $f(x)\equiv c\not\equiv 0\pmod{p}$ for any integer
$x$. Thus $S(f,p^e)$ is empty in this case.
If $p|a$ and $p\nmid b$, then there exists a unique integer
$x_0\in [1, p]$ such that $f(x_0)\equiv bx_0+c\equiv 0\pmod
p$. On the other hand, we have $f'(x_0)=2ax_0+b\equiv b\not\equiv
0\pmod p$. Then by Hensel's lemma (see, for example, \cite{[K]}),
there is a unique $p$-adic integer $s_p$ such that $f(s_p)=0$
and $s_p\equiv x_0\pmod p$. Therefore $\langle s_p\rangle_{p^e}$
is the unique solution of $f(x)\equiv0\pmod {p^e}$ in
the interval $[1, p^e]$ satisfying $\langle s_p\rangle_{p^e}\equiv
s_p\pmod{p^e}$. This completes the proof of Lemma 3.1.
\end{proof}

\noindent{\bf Lemma 3.2.} {\it Let $a$ be an odd number and let $a^{-1}$
be the inverse of $a$ in the ring $\mathbb{Z}_2$ of 2-adic integers.
For any positive integer $e$, each of the following results is true.

{\rm (i).} If either $e=2\lfloor \frac{\nu _2(D)}{2}\rfloor-1$ with
$D_4\equiv2\pmod4$ or $e\le 2\lfloor \frac{\nu _2(D)}{2}\rfloor-2$,
then
$$
S(f,2^e)= \Big\{ \Big\langle -\frac{a^{-1}b}{2}\Big\rangle_{2^{\lceil e/2\rceil}}
+m2^{\lceil e/2\rceil}: 0\le m< 2^{\lfloor
 e/2\rfloor}\Big\}.
$$

{\rm (ii).} If either $e=2\lfloor \frac{\nu _2(D)}{2}\rfloor-1$ with
$D_4\not\equiv2\pmod4$, or $e=2\lfloor \frac{\nu _2(D)}{2}\rfloor$ with
$D_4\equiv1\pmod4$, then
$$
 S(f, 2^e)=\Big\{
\Big\langle
a^{-1}\big(2^{\frac{\nu _2(D)}{2}-1}-\frac{b}{2}\big)\Big\rangle_{2^{\frac{\nu _2(D)}{2}}}
+m2^{\frac{\nu _2(D)}{2}}: 0\le m< 2^{\lfloor  e/2\rfloor}\Big\}.
$$

{\rm (iii).} If either $e=2\lfloor \frac{\nu _2(D)}{2}\rfloor$ with
$D_4\not\equiv1\pmod4$, or $e>2\lfloor \frac{\nu _2(D)}{2}\rfloor$ with
$D_4\not\equiv1\pmod8$, then $S(f, 2^e)$ is empty.

{\rm (iv).} If $D_4\equiv1\pmod8$ and $e>2\lfloor \frac{\nu _2(D)}{2}\rfloor=\nu _2(D)$, then
\begin{align*}
S(f,2^e)={\left\{\begin{array}{rl} \{\langle x_{21}\rangle_{2^e},
\langle x_{22}\rangle_{2^e}\}, & \mbox{if}\ \nu _2(D)=0,\\
 \{\langle a^{-1}(\pm X_{2^e}-\frac{b}{2})\rangle_{2^{e-\frac{\nu _2(D)}{2}}}+
m2^{e-\frac{\nu _2(D)}{2}}: 0\le m<2^{\frac{\nu _2(D)}{2}}\}, & \mbox{otherwise,}\\
 \end{array}
\right.}
\end{align*}
where $x_{21}$ and $x_{22}$ are the only two solutions of $f(x)=0$
in the ring $\mathbb{Z}_2$ of $2$-adic integers, $X_{2^e}$ denotes
the smallest root of the congruence $x^2\equiv \frac{D}{4}\pmod
{2^e}$ in the interval $[1, 2^{e-\frac{\nu _2(D)}{2}}]$.}
\begin{proof}
First, one can easily deduce from $D=b^2-4ac$ that $\nu _2(D)=0$ if $b$
is odd, and $\nu _2(D)\ge 2$ if $b$ is even. So for Cases (i) and (ii),
since $e\ge 1$, one has $\nu _2(D)>0$ and thus $b$ should be even.
If $a$ is odd and $b$ is even, then the congruence $f(x)\equiv0\pmod{2^e}$
is equivalent to
$\big(ax+\frac{b}{2}\big)^2\equiv\frac{ b^2-4ac}{4}\equiv \frac{D}{4}\pmod{2^e}.$

(i). Let  $e=2\lfloor
\frac{\nu _2(D)}{2}\rfloor-1$ with $D_4\equiv2\pmod4$ or
$e\le 2\lfloor \frac{\nu _2(D)}{2}\rfloor-2$. Then
$\frac{D}{4}\equiv 0\pmod {2^e}$. So
$y^2\equiv\frac{D}{4}\pmod{2^e}$ has exactly $2^{\lfloor
 e/2\rfloor}$ solutions: $m\cdot2^{\lceil
 e/2\rceil}$, where $1\le m \le 2^{\lfloor
 e/2\rfloor}$. Hence we can derive from
$(ax+\frac{b}{2})^2\equiv\frac{D}{4}\pmod{2^e}$ that
$ax+\frac{b}{2}\equiv m2^{\lceil  e/2\rceil}\pmod{2^e}$,
which implies that
$
x\equiv a^{-1}(m2^{\lceil e/2\rceil}-\frac{b}{2})
\equiv -\frac{a^{-1}b}{2}
+a^{-1}m2^{\lceil e/2\rceil}\pmod{2^e}
$
with $0\le m<2^{\lfloor  e/2\rfloor}$. Since $a^{-1}m$ runs
over a complete residue system modulo $2^{\lfloor e/2\rfloor}$ as $m$
does, we get
$
x\equiv \Big\langle-\frac{a^{-1}b}{2}\Big\rangle_{2^{\lceil e/2\rceil}}
+m2^{\lceil e/2\rceil}\pmod{2^e}
$
with $0\le m<2^{\lfloor  e/2\rfloor}$. Moreover, for any two integers
$m_1$ and $m_2$ satisfying $0\le m_1\ne m_2<2^{\lfloor e/2\rfloor}$, we have
$$
\Big\langle-\frac{a^{-1}b}{2}\Big\rangle_{2^{\lceil e/2\rceil}}
+m_12^{\lceil e/2\rceil}\not\equiv
\Big\langle-\frac{a^{-1}b}{2}\Big\rangle_{2^{\lceil e/2\rceil}}
+m_22^{\lceil e/2\rceil}\pmod{2^e}.
$$
So we arrive at
the desired result. Thus Part (i) is proved.

(ii). Let $e=2\lfloor \frac{\nu _2(D)}{2}\rfloor-1$ with
$D_4\not\equiv2\pmod4$ or $e=2\lfloor \frac{\nu _2(D)}{2}\rfloor$ with
$D_4\equiv1\pmod4$. Then $\nu _2(D)$ is even. From
$(ax+\frac{b}{2})^2\equiv\frac{D}{4}\pmod{2^e}$ we deduce that
$ax+\frac{b}{2}\equiv (2m+1)2^{\frac{\nu _2(D)}{2}-1}\pmod {2^e}$
with $0\le m<2^{\lfloor e/2\rfloor}$. Thus
$$x\equiv a^{-1}\Big((2m+1)2^{\frac{\nu _2(D)}{2}-1}-\frac{b}{2}\Big)\equiv
a^{-1}\big(2^{\frac{\nu _2(D)}{2}-1}-\frac{b}{2}\big)
+a^{-1}m2^{\frac{\nu _2(D)}{2}}\pmod {2^e}$$
for $0\le m< 2^{\lfloor e/2\rfloor}$. Similarly as in (i), we get
$
x\equiv \Big\langle
a^{-1}\big(2^{\frac{\nu _2(D)}{2}-1}-\frac{b}{2}\big)\Big\rangle_{2^{\frac{\nu _2(D)}{2}}}
+m2^{\frac{\nu _2(D)}{2}} \pmod{2^e}
$
with $0\le m< 2^{\lfloor e/2\rfloor}$.
On the other hand,
$$
\Big\langle
a^{-1}\big(2^{\frac{\nu _2(D)}{2}-1}-\frac{b}{2}\big)\Big\rangle_{2^{\frac{\nu _2(D)}{2}}}
+m_12^{\frac{\nu _2(D)}{2}}\not\equiv \Big\langle
a^{-1}\big(2^{\frac{\nu _2(D)}{2}-1}-\frac{b}{2}\big)\Big\rangle_{2^{\frac{\nu _2(D)}{2}}}
+m_22^{\frac{\nu _2(D)}{2}}\pmod{2^e}.
$$
for any two integers $0\le m_1\ne m_2< 2^{\lfloor  e/2\rfloor}$.
Thus
the required result
follows. So part (ii) is proved.

(iii). Let $e=2\lfloor \frac{\nu _2(D)}{2}\rfloor$ with
$D_4\not\equiv1\pmod4$ or $e>2\lfloor \frac{\nu _2(D)}{2}\rfloor$ with
$D_4\not\equiv1\pmod8$. If $\nu _2(D)=0$, then we can derive from
$D_4=D=b^2-4ac\not\equiv1\pmod8$ that $b$ and $c$ are both odd numbers. Thus
for any positive integer $n$, $f(n)\not\equiv0\pmod2$. This infers that $S(f, 2^e)$ is empty.

If $\nu _2(D)\ge 2$, since  $e=2\lfloor
\frac{\nu _2(D)}{2}\rfloor$ with $D_4\not\equiv1\pmod4$ or $e>2\lfloor
\frac{\nu _2(D)}{2}\rfloor$ with $D_4\not\equiv1\pmod8$, then
$y^2\equiv D_4\pmod{2^{e-2\lfloor \frac{\nu _2(D)}{2}\rfloor+2}}$ has no solution.
Hence there is no integer $y$ satisfying $y^2\equiv \frac{D}{4}\pmod{2^e}$.
Thus $(ax+\frac{b}{2})^2\equiv \frac{D}{4}\pmod{2^e}$ has no solution,
which means that $S(f, 2^e)$ is empty. This concludes part (iii).

(iv). Let $e>2\lfloor \frac{\nu _2(D)}{2}\rfloor$ and $D_4\equiv1\pmod8$.
If $\nu _2(D)=0$, then it follows from $D_4\equiv1\pmod8$ that $b$ is
odd and $c$ is even. Thus one has that $f(0)\equiv0\pmod 2$ and
$f(1)\equiv0\pmod 2$. On the other hand, $f'(0)\equiv
f'(1)\equiv b\not\equiv0\pmod{2}$. So by Hensel's lemma, there are
exactly two $2$-adic integers $x_{21}$ and $x_{22}$ such that
$x_{21}\equiv0\pmod2$, $x_{22}\equiv1\pmod2$ and
$f(x_{21})=f(x_{22})=0$. Thus $\langle x_{21}\rangle_{2^e}$ and
$\langle x_{22}\rangle_{2^e}$ are exactly two solutions of the
congruence $f(x)\equiv0\pmod {2^e}$ in the interval $[1, 2^e]$.

If $\nu _2(D)\ge 2$, then  $b$ is even. By the definition of $D_4$, we
know that $\nu _2(D)$ is even if $D_4\equiv 1\pmod 8$.
Since $D_4\equiv1\pmod8$ and
$e>2\lfloor \frac{\nu _2(D)}{2}\rfloor=\nu _2(D)$, it is known (see
Theorem 5.1 of page 44 in \cite{[Hu]}) that $y^2\equiv
D_4\pmod{2^{e+2-\nu _2(D)}}$ has just four solutions in the interval $[1,
2^{e+2-\nu _2(D)})$. Let $y_1$ denote the smallest solution in the interval
$[1, 2^{e+2-\nu _2(D)})$ of $y^2\equiv D_4\pmod {2^{e+2-\nu _2(D)}}$.
Evidently,  $y_1$ is odd and $y_1\in [1, 2^{e-\nu _2(D)})$. Then the
four solutions of $y^2\equiv D_4\pmod {2^{e+2-\nu _2(D)}}$ are as follows:
$y_1, 2^{e+1-\nu _2(D)}-y_1, y_1+2^{e+1-\nu _2(D)}, 2^{e+2-\nu _2(D)}-y_1.$
Thus the congruence
$y^2\equiv \frac{D}{4}\equiv 2^{\nu _2(D)-2}D_4\pmod {2^e}$
has the following solutions:
$y=2^{\frac{\nu _2(D)}{2}-1}(y_1+m_12^{e-\nu _2(D)+1})
=2^{\frac{\nu _2(D)}{2}-1}y_1+m_12^{e-\frac{\nu _2(D)}{2}}$
or
$y=2^{\frac{\nu _2(D)}{2}-1}(2^{e-\nu _2(D)+1}-y_1+m_22^{e-\nu _2(D)+1})
=2^{e-\frac{\nu _2(D)}{2}}-2^{\frac{\nu _2(D)}{2}-1}y_1+m_22^{e-\frac{\nu _2(D)}{2}},$
where $0\le m_1, m_2<2^{\frac{\nu _2(D)}{2}}$ are integers. Now let
$X_{2^e}=2^{\frac{\nu _2(D)}{2}-1}y_1$. Then from
$(ax+\frac{b}{2})^2\equiv \frac{D}{4}\equiv 2^{\nu _2(D)-2}D_4\pmod{2^e}$
we get that
$ax+\frac{b}{2}\equiv \pm X_{2^e}+m2^{e-\frac{\nu _2(D)}{2}}\pmod{2^e},$
which implies that
$x\equiv\langle a^{-1}(\pm X_{2^e}-\frac{b}{2})
\rangle_{2^{e-\frac{\nu _2(D)}{2}}}+a^{-1}m2^{e-\frac{\nu _2(D)}{2}} \pmod {2^e}$
for any integer $m$ with $0\le m<2^{\frac{\nu _2(D)}{2}}$.
Since $a^{-1}m$ runs over a complete residue system
modulo $2^{\frac{\nu _2(D)}{2}}$ as $m$
does so, we obtain that
$
x\equiv\langle a^{-1}(\pm X_{2^e}-\frac{b}{2})
\rangle_{2^{e-\frac{\nu _2(D)}{2}}}+m2^{e-\frac{\nu _2(D)}{2}} \pmod {2^e}
$
for any integer $m$ with $0\le m<2^{\frac{\nu _2(D)}{2}}$.
One can easily check that all the $2^{\frac{\nu _2(D)}{2}+1}$ elements of
the set $\{\langle a^{-1}(\pm X_{2^e}
-\frac{b}{2})\rangle_{2^{e-\frac{\nu _2(D)}{2}}}+m2^{e-\frac{\nu _2(D)}{2}}:
0\le m<2^{\frac{\nu _2(D)}{2}}\}$ are pairwise incongruent modulo $2^e$. Thus the
desired result follows immediately. The proof of Lemma 3.2 is complete.
\end{proof}

\noindent {\bf Lemma 3.3.} {\it Let $p$ be an odd prime with $p\nmid
a$, and let $D_p$ be defined as in (1.3). By $(2a)^{-1}$ we denote the inverse
of $2a$ in the ring $\mathbb{Z}_p$ of $p$-adic integers.
For any positive integer $e$, each of the following results is true.

{\rm (i).} If $e\le \nu _p(D)$, then
$S(f, p^e)=\big\{\langle -(2a)^{-1}b\rangle_{p^{\lceil e/2\rceil}}+mp^{\lceil
e/2\rceil}: 0\le m< p^{\lfloor e/2\rfloor}\big\}.$

{\rm (ii).} If $e>\nu _p(D)$ with either $\nu _p(D)$ being odd or $
(\frac{D_p}{p})=-1$, then $S(f, p^e)$ is empty.

{\rm (iii).} If $e>\nu _p(D)$ with $\nu _p(D)$ being even and $(\frac{D_p}{p})=1$, then
$$
S(f, p^e)=\Big\{\Big\langle (2a)^{-1}\big(
\pm X_{p^e}-b\big)\Big\rangle_{p^{e-\frac{\nu _p(D)}{2}}}+
mp^{e-\frac{\nu _p(D)}{2}}: 0\le
m<p^{\frac{\nu _p(D)}{2}}\Big\},
$$
where $X_{p^{e}}$ is the smallest solution of $x^2\equiv
D\pmod{p^{e}}$ in the interval $[1, p^{e-\frac{\nu _p(D)}{2}}]$. }
\begin{proof} Since $p\nmid 2a$, $f(x)\equiv0\pmod{p^e}$ is
equivalent to $(2ax+b)^2\equiv D\pmod{p^e}.$

{\rm (i).} Let $e\le \nu _p(D)$. Then one has $2ax+b\equiv mp^{\lceil e/2\rceil}\pmod{p^e}$
for some integers $1\le m\le p^{\lfloor e/2\rfloor}$. Hence
$x\equiv (2a)^{-1}(mp^{\lceil e/2\rceil}-b)\pmod{p^e}$ for
$1\le m\le p^{\lfloor e/2\rfloor}$.
Moreover, we have $(2a)^{-1}(m_1p^{\lceil e/2\rceil}-b)\not\equiv
(2a)^{-1}(m_2p^{\lceil e/2\rceil}-b)\pmod{p^e}$
for any two integers $m_1$ and $m_2$ with $1\le m_1\ne m_2\le
p^{\lfloor e/2\rfloor}$. On the other hand, since $(2a)^{-1}m$ runs over
a complete residue system modulo $p^{\lfloor e/2\rfloor}$ as $m$ does, we get that
$x\equiv \langle -(2a)^{-1}b\rangle_{p^{\lceil e/2\rceil}}+mp^{\lceil
e/2\rceil}\pmod{p^e}$ with $ 0\le m< p^{\lfloor e/2\rfloor}$.
Thus we derive the required result immediately.

{\rm (ii).} Let $e>\nu _p(D)$ with either $\nu _p(D)$ odd or $
(\frac{D_p}{p})=-1$. Suppose that there is an integer $n_0$ such that
$n_0^2\equiv D\pmod {p^e}$. Then $n_0^2\equiv p^{\nu _p(D)}D_p\pmod {p^e}$.
Since $e>\nu _p(D)$, $p^{\lceil \frac{\nu _p(D)}{2}\rceil}$ divides $n_0$.
If $\nu _p(D)$ is odd, then
$n_0^2\equiv p^{\nu _p(D)+1}\frac{n_0^2}{p^{\nu _p(D)+1}}\equiv p^{\nu _p(D)}D_p\pmod{p^e}.$
Hence $p\cdot\frac{ n_0^2}{p^{\nu _p(D)+1}}\equiv D_p\pmod{p^{e-\nu _p(D)}},$
which is a contradiction. If $\nu _p(D)$ is even and $(\frac{D_p}{p})=-1$, then
$\Big(\frac{n_0}{p^{\nu _p(D)/2}}\Big)^2\equiv D_p\pmod {p^{e-\nu _p(D)}},$
which is impossible since $(\frac{D_p}{p})=-1$.
Thus there is no integer $y$ such that $y^2\equiv D\pmod{p^e}$.
It follows immediately that the congruence $(2ax+b)^2\equiv D\pmod{p^e}$
has no solution and $S(f,p^e)$ is empty.

(iii). Let $e>\nu _p(D)$ with $\nu _p(D)$ even and $(\frac{D_p}{p})=1$.
Then by Hensel's lemma, the congruence $y^2\equiv D_p\pmod {p^{e-\nu _p(D)}}$
has exactly two solutions $y_0$ and $p^{e-\nu _p(D)}-y_0$
in the interval $[1, p^{e-\nu _p(D)}]$. Thus for all integers
$0\le m<p^{\frac{\nu _p(D)}{2}}$, we have
$y=p^{\frac{\nu _p(D)}{2}}(y_0+mp^{e-\nu _p(D)})
=p^{\frac{\nu _p(D)}{2}}y_0+mp^{e-\frac{\nu _p(D)}{2}}$
and
$y=p^{e-\frac{\nu _p(D)}{2}}-p^{\frac{\nu _p(D)}{2}}y_0+mp^{e-\frac{\nu _p(D)}{2}},$
are solutions of the congruence $y^2\equiv D\pmod {p^e}$.
Let now $X_{p^e}=p^{\frac{\nu _p(D)}{2}}y_0$. Then we derive that
$2ax+b\equiv \langle \pm X_{p^e}\rangle_{p^{e-\frac{\nu _p(D)}{2}}}+
mp^{e-\frac{\nu _p(D)}{2}}\pmod {p^e}.$
Since $(2a)^{-1}m$ runs over the complete residue system as $m$ does, we get
$x\equiv \Big\langle (2a)^{-1}\big(
\pm X_{p^e}-b\big)\Big\rangle_{p^{e-\frac{\nu _p(D)}{2}}}+
mp^{e-\frac{\nu _p(D)}{2}} \pmod{p^e}
$
where $0\le m<p^{\frac{\nu _p(D)}{2}}$. Obviously, any two classes of the above
$2p^{\frac{\nu _p(D)}{2}}$ residue classes are distinct modulo $p^e$.
This concludes the desired result. Lemma 3.3 is proved.
\end{proof}

Once we determine  the set of solutions of $f(x)\equiv 0\pmod{p^e}$,
we naturally want to know more information about these
solutions. First, we introduce the following concepts,
which are important ingredients in the process of determining the local periods.

\noindent {\bf Definition 3.4.} Let $e$ be a nonnegative integer.
If $S(f, p^e)$ is nonempty, then for any $x_1, x_2\in S(f, p^e)$, we define
the {\it distance}  of $x_1$ and $x_2$ by
$$d_{p^e}(x_1, x_2):=\min\{\langle x_1-x_2\rangle_{p^e}, \langle x_2-x_1\rangle_{p^e}\}.$$

Clearly, for any $x_1, x_2\in S(f, p^e)$, $d_{p^e}(x_1, x_2)$ equals
$\min\{|x_1-x_2|, p^e-|x_1-x_2|\}$ if $x_1\ne x_2$, and is $p^e$ if $x_1=x_2$.

\noindent {\bf Definition 3.5.} Let $e$ be a nonnegative integer.
We define {\it the minimal distance}, denoted by $d_{p^e}$, among
the solutions of $f(x)\equiv0\pmod{p^e}$ as follows: $d_{p^0}:=1$,
and for $e\ge 1$,
\begin{align*}
d_{p^e}:={\left\{
\begin{array}{rl}
\min\{d_{p^e}(x_i, x_j): x_i, x_j\in S(f, p^e)\}, & \mbox{if}\
S(f, p^e)\  \mbox{is\ nonempty},\\
\infty, & \mbox{if}\  S(f, p^e)\  \mbox{is\ empty}.
\end{array}
\right.}
\end{align*}

In what follows, we study the arithmetic properties of the minimal distance $d_{p^e}$.

\noindent {\bf Lemma 3.6.} {\it Let $S(f, p^e)$ be nonempty.
Then there exists a positive integer $n$ such that
$\nu _p(f(n))\ge e$ and $\nu _p(f(n+d_{p^e}))\ge e$.
Further, if $d_{p^e}<d_{p^{e+1}}$, then there
is a positive integer $m$ such that $\nu _p(f(m))=e$.}
\begin{proof} First let $|S(f, p^e)|=1$. Then $S(f, p^e)$ contains only one element,
saying $x_0$, and $d_{p^e}=p^e$. Thus one can pick $n=x_0$ to arrive at the desired result.
Lemma 3.6 is true in this case. Now let $|S(f, p^e)|\ge 2$. Then by Definitions 3.5 and 3.4,
there are $x_1, x_2\in S(f, p^e)$ with $x_1<x_2$ such that $d_{p^e}=d_{p^e}(x_1, x_2)$.
It follows that
\begin{align}\label{2.3}
d_{p^e}=\min(x_2-x_1, p^e+x_1-x_2).
\end{align}

If $x_2-x_1\le \frac{1}{2}p^{e}$, by (\ref{2.3}) we have $d_{p^e}=x_2-x_1$.
Take $n=x_1$. Then we have $\nu _p(f(n))=\nu _p(f(x_1))\ge e$ and $\nu _p(f(n+d_{p^e}))=\nu _p(f(x_2))\ge e$
as required.

If $x_2-x_1>\frac{1}{2}p^{e}$, then by (\ref{2.3}) we have $d_{p^e}=p^e+x_1-x_2$. Hence
taking $n=x_2$ gives that $n+d_{p^e}=x_1+p^e$. But $f(x_1+p^e)\equiv f(x_1)\pmod {p^e}$
and $f(x_1)\equiv 0\pmod {p^e}$. So we have $f(n+d_{p^e})\equiv 0\pmod {p^e}$.
The first part of Lemma 3.6 is proved.

Now suppose that $d_{p^e}<d_{p^{e+1}}$.
Since $S(f, p^e)$ is nonempty, by the first part
we can find a positive integer $n_1$ such that
$\nu _p(f(n_1))\ge e$ and $\nu _p(f(n_1+d_{p^e}))\ge e$. Suppose that
$\nu _p(f(n_1))>e$ and $\nu _p(f(n_1+d_{p^e}))>e$. Then
$\nu _p(f(n_1))\ge e+1$ and $\nu _p(f(n_1+d_{p^e}))\ge e+1$.
Thus $\nu _p(f(\langle n_1\rangle _{p^{e+1}}))\ge e+1$ and
$\nu _p(f(\langle n_1+d_{p^e}\rangle _{p^{e+1}}))\ge e+1$.
That is, $\langle n_1\rangle _{p^{e+1}}$ and
$\langle n_1+d_{p^e}\rangle _{p^{e+1}}$ are belonging to
the set $S(f, p^{e+1})$. Since $\langle n_1+d_{p^e}\rangle _{p^{e+1}}-
\langle n_1\rangle _{p^{e+1}}\equiv d_{p^e}\pmod {p^{e+1}}$, we get that
$|\langle n_1+d_{p^e}\rangle _{p^{e+1}}-
\langle n_1\rangle _{p^{e+1}}|=d_{p^e}$ or $p^{e+1}-d_{p^e}$.
It then follows that
$d_{p^{e+1}}\le d_{p^{e+1}}(\langle n_1\rangle _{p^{e+1}},
\langle n_1+d_{p^e}\rangle _{p^{e+1}})
=\min (d_{p^e}, p^{e+1}-d_{p^e})\le d_{p^e}$
which contradicts with the assumption $d_{p^e}<d_{p^{e+1}}$.
Therefore we have either $\nu _p(f(n_1))=e$ or $\nu _p(f(n_1+d_{p^e}))=e$.
This concludes Lemma 3.6.
\end{proof}

\noindent {\bf Lemma 3.7.} {\it For any given prime number $p$,
the sequence $\{d_{p^e}\}_{e=0}^\infty $ is nondecreasing.}
\begin{proof} Evidently, it is enough to prove $d_{p^{e+1}}
\ge d_{p^e}$ for any nonnegative integer $e$. Now let $e$ be any given nonnegative
integer. If $S(f, p^{e+1})$ is empty, then $d_{p^{e+1}}=\infty$
and so $d_{p^{e+1}}\ge d_{p^e}$ as desired. In what follows we assume that
$S(f, p^{e+1})$ is not empty.

Clearly, it suffices to prove that $d_{p^{e+1}}(x_1, x_2)\ge d_{p^e}$
for any two elements $x_1, x_2\in S(f, p^{e+1})$, from which we derive that
$d_{p^{e+1}}\ge d_{p^e}$. Taking any two elements $x_1, x_2\in S(f, p^{e+1})$,
then we have that $x_1=y_1+t_1p^e$ and $x_2=y_2+t_2p^e$ for some integers
$t_1, t_2, y_1, y_2$ with $0\le t_1, t_2\le p-1$ and $y_1, y_2\in S(f, p^e)$.
It is clear that
$$d_{p^{e+1}}(x_1, x_2)=\min\{ \langle(t_1-t_2)p^e\rangle_{p^{e+1}},
\langle(t_2-t_1)p^e\rangle_{p^{e+1}}\}\ge p^e=d_{p^e}(y_1, y_2)\ge d_{p^e}$$
if $y_1=y_2$. Let now $y_1\ne y_2$. Then $x_1\ne x_2$. We have
$d_{p^e}(y_1, y_2)=\min\{ |y_1-y_2|, p^e-|y_1-y_2|\}$ and
$d_{p^{e+1}}(x_1, x_2)=\min\{|y_1-y_2+(t_1-t_2)p^e|,
p^{e+1}-|y_1-y_2+(t_1-t_2)p^e|\}.$

If $t_1=t_2$, then $d_{p^{e+1}}(x_1, x_2)=|y_1-y_2|\ge d_{p^e}(y_1, y_2)\ge d_{p^e}$.
If $t_1\ne t_2$, then
$$|y_1-y_2+(t_1-t_2)p^e|\ge |(t_1-t_2)p^e|-|y_1-y_2|
\ge p^e-|y_1-y_2|\ge d_{p^e}(y_1, y_2)\ge d_{p^e}$$
and
\begin{align*}
p^{e+1}-|y_1-y_2+(t_1-t_2)p^e|
&\ge p^{e+1}-|y_1-y_2|-|(t_1-t_2)p^e|\\
&\ge p^e-|y_1-y_2|\ge d_{p^e}(y_1, y_2)\ge d_{p^e}.
\end{align*}
It follows that $d_{p^{e+1}}(x_1, x_2)\ge d_{p^e}$.
So Lemma 3.7 is proved.
\end{proof}

\noindent{\bf Lemma 3.8.} {\it  If $p|a$,  then for any positive
integer $e$, we have
$
d_{p^e}={\left\{
  \begin{array}{rl}
\infty, & \mbox{if}\ p|b,\\
p^e, & \mbox{if}\ p\nmid b.
 \end{array}
\right.}
$}

\begin{proof}
It follows immediately from Lemma 3.1 that Lemma 3.8 is true.
\end{proof}

\noindent{\bf Lemma 3.9.} {\it Let $a$ be odd and $e$ be a positive integer.

{\rm (i).} If $e=\nu _2(D)$ with $D_4\equiv1\pmod4$ or $e\le 2\lfloor
\frac{\nu _2(D)}{2}\rfloor-1$, then $d_{2^e}=2^{\lceil  e/2\rceil}$. Also
$d_{2^e}$ equals the smallest positive root of $a^2x^2-D\equiv0\pmod{2^e}$.
Moreover, the distance between any two
distinct solutions of $f(x)\equiv0\pmod {2^e}$ is divisible by
$2^{\lceil  e/2\rceil}$ if $e\ge 2$.

{\rm (ii).} If $e= 2\lfloor \frac{\nu _2(D)}{2}\rfloor$ with
$D_4\not\equiv1\pmod4$  or $e>2\lfloor \frac{\nu _2(D)}{2}\rfloor$ with
$D_4\not\equiv1\pmod8$, then $d_{2^e}=\infty.$

{\rm (iii).} If $e>2\lfloor \frac{\nu _2(D)}{2}\rfloor=\nu _2(D)$ with
$D_4\equiv1\pmod8$, then $d_{2^{e}}$ is equal to the smallest
positive root of the congruence $a^2x^2-D\equiv0\pmod {2^{e+1}}$.}
\begin{proof}
(i). It is easy to check that $d_{2^e}=2^e=2^{\lceil e/2\rceil}$ if $e=1$.
So it is enough to prove part (i) for the case $e\ge 2$.
In what follows we let $e\ge 2$.

If $e=2\lfloor \frac{\nu _2(D)}{2}\rfloor-1$ with
$D_4\equiv2\pmod4$ or $e\le 2\lfloor \frac{\nu _2(D)}{2}\rfloor-2$,
then by Lemma 3.2 (i), $\langle -\frac{a^{-1}b}{2}\rangle_{2^{\lceil e/2\rceil}}
+i2^{\lceil e/2\rceil}$ and $\langle -\frac{a^{-1}b}{2}\rangle_{2^{\lceil e/2\rceil}}
+2^{\lceil e/2\rceil}$ are two roots
of $f(x)\equiv0\pmod {2^e}$, where
$ 0\le i\ne j<2^{\lfloor e/2\rfloor}$. It is easy to see that
\begin{align}\label{2.4}
\langle -\frac{a^{-1}b}{2}\rangle_{2^{\lceil e/2\rceil}}
+i2^{\lceil e/2\rceil}
-(\langle -\frac{a^{-1}b}{2}\rangle_{2^{\lceil e/2\rceil}}
+j2^{\lceil e/2\rceil})
\equiv (i-j)2^{\lceil  e/2\rceil} \pmod{2^e}.
\end{align}
We can find two integers $i_0$ and
$j_0$ with $0\le i_0\ne j_0< 2^{\lfloor  e/2\rfloor}$ such
that $(i_0-j_0)=1$. Then by (\ref{2.4}), we have
$d_{2^e}=2^{\lceil e/2\rceil}$ as required.

If either $e=2\lfloor \frac{\nu _2(D)}{2}\rfloor-1$ with
$D_4\not\equiv2\pmod4$ or $e=2\lfloor \frac{\nu _2(D)}{2}\rfloor =\nu _2(D)$ with
$D_4\equiv1\pmod4$, then $\nu _2(D)\ge 2$ is even in this case. For any
two integers $i$ and $j$ with $0\le i\ne j< 2^{\frac{\nu _2(D)}{2}}$, we
have
\begin{align}\label{2.5}
\Big\langle
a^{-1}\big(2^{\frac{\nu _2(D)}{2}-1}-\frac{b}{2}\big)\Big\rangle_{2^{\frac{\nu _2(D)}{2}}}
+i2^{\frac{\nu _2(D)}{2}}-\Big\langle a^{-1}\big(2^{\frac{\nu _2(D)}{2}-1}-\frac{b}{2}
\big)\Big\rangle_{2^{\frac{\nu _2(D)}{2}}}-j2^{\frac{\nu _2(D)}{2}}=(i-j)2^{\frac{\nu _2(D)}{2}}.
\end{align}
Then by Lemma 3.2 (ii), we have $d_{2^e}=2^{\frac{\nu _2(D)}{2}}=2^{\lceil
 e/2\rceil}$ as desired. Evidently, under the assumptions of part (i), we have
$e\le \nu _2(D)$. So $2^{\lceil  e/2\rceil}$ is the smallest positive
root of $a^2x^2-D\equiv0\pmod{2^e}$. Finally, by (\ref{2.4}) and (\ref{2.5}),
$2^{\lceil  e/2\rceil}$ divides the distance between
any two distinct solutions of $f(x)\equiv0\pmod {2^e}$.

(ii). By Lemma 3.2 (iii), $S(f, 2^e)$ is empty in this case, which means
that $d_{2^e}=\infty$.

(iii). Since $D_4\equiv 1\pmod 8$, $\nu _2(D)$ is even. If $\nu _2(D)=0$,
then $D=D_4$. Clearly $d_2=1$ is the
smallest positive root of $a^2x^2-D\equiv0\pmod{2^2}$. In what
follows, we only need to consider the case either $\nu _2(D)=0$ and $e\ge 2$, or
$\nu _2(D)\ge 2$ and $e\ge \nu _2(D)+1$.

For any integer $x_1$ satisfying $f(x_1)\equiv0\pmod {2^{e}}$, we have
$(2ax_1+b)^2\equiv D\pmod {2^{e+2}}.$
Hence $(2ax_1+b)^2\equiv D\pmod{2^{e+1}}$. Note that the discriminant of $x^2-D$
equals $4D$. Since
$\nu _2(4D)=\nu _2(D)+2\ {\rm and} \ \frac{4D}{2^{\nu _2(D)+2}}=D_4\equiv1\pmod8,$
then by Lemma 3.2 (iv), we obtain the exactly $2^{\frac{\nu _2(D)+2}{2}+1}=2^{\frac{\nu _2(D)}{2}+2}$
roots of $x^2-D\equiv0\pmod {2^{e+1}}$. Now let $y_{e+1}$ be the smallest
positive root of $a^2x^2-D\equiv0\pmod {2^{e+1}}$. Then $y_{e+1}\in [1, 2^{e-1-\nu _2(D)/2})$
and $2^{\frac{\nu _2(D)}{2}}$ divides $y_{e+1}$. We claim that
the set of solutions  of $x^2-D\equiv0\pmod {2^{e+1}}$ in the interval $[1, 2^{e+1}]$
is given as follows:
\begin{align}\label{2.7}
\{\langle \pm ay_{e+1}\rangle_{2^{e-\frac{\nu _2(D)}{2}}}
+m2^{e-\frac{\nu _2(D)}{2}}: 0\le m<2^{\frac{\nu _2(D)}{2}+1}\}.
\end{align}
On the one hand, since $2^{\frac{\nu _2(D)}{2}}$ divides $y_{e+1}$ and $e>\nu _2(D)$,
one can easily check that each one in the set (\ref{2.7}) satisfies
$x^2-D\equiv0\pmod {2^{e+1}}$. On the other hand,
it is clear that any two elements
in the set (\ref{2.7}) are incongruent modulo $2^{e+1}$ and the set (\ref{2.7}) holds
$2^{\frac{\nu _2(D)}{2}+2}$ elements. Thus (\ref{2.7}) gives all the solutions
of the congruence $x^2-D\equiv0\pmod {2^{e+1}}$ in the
interval $[1, 2^{e+1}]$. The claim is proved.

Now from the claim we get that
$2ax_1+b\equiv \pm ay_{e+1}+m_02^{e-\frac{\nu _2(D)}{2}} \pmod {2^{e+1}}$
for some $0\le m_0<2^{\frac{\nu _2(D)}{2}+1}$, which implies that
$2ax_1+b\equiv  ay_{e+1}\  \mbox{or} \ -ay_{e+1}\pmod{2^{e-\frac{\nu _2(D)}{2}}}.$

In what follows we show that $d_{2^e}=y_{e+1}$. Since the proof for the case
$2ax_1+b\equiv -ay_{e+1}\pmod{2^{e-\frac{\nu _2(D)}{2}}}$ is similar as that of the
case $2ax_1+b\equiv ay_{e+1}\pmod{2^{e-\frac{\nu _2(D)}{2}}}$, we only treat with
the latter case. Let $2ax_1+b\equiv ay_{e+1}\pmod{2^{e-\frac{\nu _2(D)}{2}}}$.
Then we have $$a(x_1-y_{e+1})^2+b(x_1-y_{e+1})+c\equiv
-(2ax_1+b)y_{e+1}+ay_{e+1}^2$$
$$\equiv -(ay_{e+1}+t2^{e-\frac{\nu _2(D)}{2}})y_{e+1}+ay_{e+1}^2\equiv
-ty_{e+1}2^{e-\frac{\nu _2(D)}{2}}\pmod {2^e}
$$
for some integer $t$. Since $2^{\frac{\nu _2(D)}{2}}\big|y_{e+1}$, we then derive that
$f(x_1-y_{e+1})\equiv-ty_{e+1}2^{e-\frac{\nu _2(D)}{2}}\equiv 0\pmod{2^{e}}.$
So $x_1-y_{e+1}$ is a solution of $f(x)\equiv 0\pmod {2^e}$.
But $x_1$ is a solution of $f(x)\equiv 0\pmod {2^e}$.
It follows that $d_{2^e}\le d_{2^e}(x_1, x_1-y_{e+1})\le y_{e+1}$.

Noticing that
$
f(x_1+m2^{e-\frac{\nu _2(D)}{2}})=f(x_1)+(2ax_1+b)m2^{e-\frac{\nu _2(D)}{2}}+am^22^{2e-\nu _2(D)},
$
$2^{\frac{\nu _2(D)}{2}}|(2ax_1+b)$ and $e>\nu _2(D)$, we deduce that
$f(x_1+m2^{e-\frac{\nu _2(D)}{2}})\equiv 0\pmod {2^e}$
for any integer $m$. Replacing $x_1$ by $x_1-y_{e+1}$ gives that
$f(x_1-y_{e+1}+m2^{e-\frac{\nu _2(D)}{2}})\equiv 0\pmod {2^e}$
for any integer $m$. Thus by Lemma 3.2 (iv), one knows that
$\langle x_1 \rangle _{2^{e-\frac{\nu _2(D)}{2}}},
\langle x_1-y_{e+1} \rangle _{2^{e-\frac{\nu _2(D)}{2}}}$
are the only two distinct solutions of $f(x)\equiv 0\pmod {2^e}$
in the interval $[1, 2^{e-\frac{\nu _2(D)}{2}}]$.

Since
$
\langle x_1 \rangle _{2^{e-\frac{\nu _2(D)}{2}}}
-\langle x_1-y_{e+1} \rangle _{2^{e-\frac{\nu _2(D)}{2}}}
\equiv y_{e+1} \pmod {2^{e-\frac{\nu _2(D)}{2}}},
$
one can easily check that
$\langle x_1 \rangle _{2^{e-\frac{\nu _2(D)}{2}}}-\langle x_1-y_{e+1}
\rangle _{2^{e-\frac{\nu _2(D)}{2}}}$
is a solution of $a^2x^2-D\equiv 0\pmod {2^{e+1}}$.
It then follows immediately that so is
\begin{align}\label{2.8}
\langle \pm \langle x_1 \rangle _{2^{e-\frac{\nu _2(D)}{2}}}\mp \langle x_1-y_{e+1}
\rangle _{2^{e-\frac{\nu _2(D)}{2}}}+l2^{e-\frac{\nu _2(D)}{2}}\rangle _{2^e}
\end{align}
for any integer $l$. However, by the definition of $d_{2^e}$ and Lemma 3.2 (iv),
$d_{2^e}$ must be of the form (\ref{2.8}) for some integer $l$.
Therefore $d_{2^e}$ is a solution of $a^2x^2-D\equiv 0\pmod {2^{e+1}}$.
Then by the minimality of $y_{e+1}$, we have $d_{2^{e}}\ge y_{e+1}$.
So we get that $d_{2^{e}}=y_{e+1}$ as desired. Part (iii) is proved.
The proof of Lemma 3.9 is complete.
\end{proof}
\noindent {\bf Lemma 3.10.} {\it Let $p$ be an odd prime with $p\nmid
a$, and let $D_p$ be defined as in (1.3).

{\rm (i).} If $e\le \nu _p(D)$, then $d_{p^e}=p^{\lceil  e/2\rceil}$. Further,
$d_{p^e}$ is equal to the smallest positive root of the congruence
$a^2x^2-D\equiv0\pmod{p^e}$. Moreover, the distance between any two
distinct solutions of the congruence $f(x)\equiv0\pmod{p^e}$ is divisible
by $p^{\lceil e/2\rceil}$ if $e\ge 2$.

{\rm (ii).} If  $e>\nu _p(D)$ with either $\nu _p(D)$ being odd or $
(\frac{D_p}{p})=-1$, then $d_{p^e}=\infty.$

{\rm (iii).} If $e>\nu _p(D)$ with $\nu _p(D)$ being even and $(\frac{D_p}{p})=1$,
then $d_{p^e}$ equals the smallest positive root of the
congruence $a^2x^2-D\equiv0\pmod{p^e}$.}
\begin{proof}
(i). Obviously, part (i) is true if $e=1$. So we only need to show
that part (i) holds for the case $e\ge 2$. Now let $e\ge 2$.
For any given two integers $m_1$ and $m_2$ with $0\le m_1\ne
m_2<p^{\lfloor e/2\rfloor}$, we have
$
\langle -(2a)^{-1}b\rangle_{p^{\lceil e/2\rceil}}+m_1p^{\lceil
e/2\rceil}-(\langle -(2a)^{-1}b\rangle_{p^{\lceil e/2\rceil}}+m_2p^{\lceil
e/2\rceil})\equiv(m_1-m_2)p^{\lceil e/2\rceil} \pmod{p^e}.
$
It is easy to derive from the above congruence and Lemma 3.3 (i) that the distance
between any two distinct solutions of $f(x)\equiv0\pmod{p^e}$ is
divisible by $p^{\lceil e/2\rceil}$.

Pick two integers $m_1'$ and $m_2'$ in the interval $[0, p^{\lfloor e/2\rfloor})$
such that $m_1'-m_2'=1$. Thus $d_{p^e}=p^{\lceil e/2\rceil}$. But $e\le \nu _p(D)$.
So $p^{\lceil e/2\rceil}$ is the smallest solution of the congruence
$a^2x^2-D\equiv0\pmod{p^e}$. Part (i) is proved.

(ii). Since $e>\nu _p(D)$ with $\nu _p(D)$ being odd or $ (\frac{D_p}{p})=-1$, by Lemma 3.3
we know that $S(f, p^e)$ is empty. Thus $d_{p^e}=\infty$ as desired.

(iii). Since $e>\nu _p(D)$ with $\nu _p(D)$ being even and $(\frac{D_p}{p})=1$,
by Lemma 3.3 (iii) the congruence $a^2x^2-D\equiv0\pmod {p^e}$ has exactly
$2p^{\frac{\nu _p(D)}{2}}$ roots in any complete residue system modulo
$p^e$. Let $y_e$ be the smallest positive root of the congruence
$a^2x^2-D\equiv0\pmod{p^e}$. Then Lemma 3.3 (iii) applied to the
congruence $a^2x^2-D\equiv0\pmod{p^e}$ gives that
$y_e\in [1, p^{e-\frac{\nu _p(D)}{2}}]$ and $\nu _p(y_e)=\nu _p(X_{p^e})=\frac{\nu _p(D)}{2}$.
One can easily check that
$\langle ay_e\rangle_{p^{e-\frac{\nu _p(D)}{2}}}$ and
$\langle -ay_e\rangle_{p^{e-\frac{\nu _p(D)}{2}}}$
are the only two solutions of $x^2-D\equiv0\pmod {p^e}$ in the interval
$[1, p^{e-\frac{\nu _p(D)}{2}}]$. So by Lemma 3.3 (iii), the following set
\begin{align}\label{2.9}
\big\{  \langle\pm ay_e\rangle_{p^{e-\frac{\nu _p(D)}{2}}}+mp^{e-
\frac{\nu _p(D)}{2}}: 0\le m<p^{\frac{\nu _p(D)}{2}}\big\}
\end{align}
is exactly the set of all the solutions of the congruence
$x^2-D\equiv0\pmod{p^e}$ in the interval $[1, p^e]$. For
any solution $x_0$ of $f(x_0)\equiv0\pmod {p^e}$, one
has $(2ax_0+b)^2\equiv  D\pmod {p^e}$. By (\ref{2.9}), we have
$2ax_0+b\equiv \langle\pm ay_e\rangle_{p^{e-\frac{\nu _p(D)}{2}}}+mp^{e-
\frac{\nu _p(D)}{2}}\pmod {p^e}$ for some $0\le m<p^{e-
\frac{\nu _p(D)}{2}}$,
which implies that $2ax_0+b\equiv ay_e\ \mbox{or}\ -ay_e\pmod{p^{e-\frac{\nu _p(D)}{2}}}$.

Now we prove that $d_{p^e}=y_e$. We only need to give the the proof for the case
$2ax_0+b\equiv ay_e\pmod{p^{e-\frac{\nu _p(D)}{2}}}$, since the proof for the case
$2ax_0+b\equiv -ay_e\pmod{2^{e-\frac{\nu _2(D)}{2}}}$ is similar. Let now
$2ax_0+b\equiv ay_e\pmod{p^{e-\frac{\nu _p(D)}{2}}}$. Then
$2ax_0+b=ay_e+mp^{e-\frac{\nu _2(D)}{2}}$ for some integer $m$.
Using the fact that $\nu _p(y_e)=\frac{\nu _p(D)}{2}$, we get
$f(x_0-y_e)\equiv
-(2ax_0+b)y_e+ay_e^2 \equiv 0\pmod{p^{e}},$
i.e., $x_0-y_e$ is a solution of $f(x)\equiv0\pmod{p^e}$.
Note that $x_0$ is also a solution of $f(x)\equiv0\pmod{p^e}$.
Therefore, $d_{p^e}\le d_{p^e}(x_0, x_0-y_e)\le y_e$.

We can check that
$\langle x_0\rangle_{p^{e-\frac{\nu _p(D)}{2}}},
\langle x_0-y_e\rangle_{p^{e-\frac{\nu _p(D)}{2}}}$
are the only two distinct
solutions of $f(x)\equiv 0\pmod {p^e}$ in the interval $[1, p^{e-\frac{\nu _p(D)}{2}}]$.
Thus by Lemma 3.3 (iii), $S(f, p^e)$ is equal to the following union:
$$ \{\langle x_0\rangle_{p^{e-\frac{\nu _p(D)}{2}}}
+mp^{e-\frac{\nu _p(D)}{2}}: 0\le m<p^{\frac{\nu _p(D)}{2}}\}\cup
\{\langle x_0-y_e\rangle_{p^{e-\frac{\nu _p(D)}{2}}}
+mp^{e-\frac{\nu _p(D)}{2}}: 0\le m<p^{\frac{\nu _p(D)}{2}}\}.$$
Then $d_{p^e}$ must be of the form
$\Big\langle \pm \langle x_0\rangle_{p^{e-\frac{\nu _p(D)}{2}}}
\mp\langle x_0-y_e\rangle_{p^{e-\frac{\nu _p(D)}{2}}}
+mp^{e-\frac{\nu _p(D)}{2}}\Big\rangle_{p^e}$
for some integer $m$. It follows from Lemma 3.3 (iii) that
$\pm y_e+tp^{e-\frac{\nu _p(D)}{2}}$ is the root of $a^2x^2-D\equiv0\pmod{p^e}$
for any integer $t$. On the other hand, since
$\Big\langle \pm \langle x_0\rangle_{p^{e-\frac{\nu _p(D)}{2}}}
\mp\langle x_0-y_e\rangle_{p^{e-\frac{\nu _p(D)}{2}}}
+mp^{e-\frac{\nu _p(D)}{2}}\Big\rangle_{p^e}
\equiv \pm y_e\pmod{p^{e-\frac{\nu _p(D)}{2}}}$
for any integer $m$, $d_{p^e}$ is a solution of $a^2x^2-D\equiv0\pmod{p^e}$.
Since $y_e$ is the smallest positive root of the congruence
$a^2x^2-D\equiv0\pmod{p^e}$, we have $d_{p^{e}}\ge y_e$. So
$d_{p^{e}}=y_e$ as required. Part (iii) is proved.
This completes the proof of Lemma 3.10.
\end{proof}

Associated to $f$, we define a subset $\mathcal {K}_f$ of the set
${\mathbb N}^*$ of positive integers by
$\mathcal {K}_f:=\{j\in {\mathbb N}^*: D\ne a^2i^2 \
{\rm for \ any \ integer}\ i \ {\rm such \ that} \ 1\le i\le j\}.$
Clearly $\mathcal {K}_f$ is empty if and only if $D=a^2$. Further,
the condition that $D\ne a^2i^2$ for all integers $i$ with
$1\le i\le k$ presented in Theorem 1.2 is equivalent to that
$\mathcal {K}_f$ is nonempty and $k\in \mathcal {K}_f$.
The following result describes the form of $\mathcal{K}_f$.

\noindent {\bf Lemma 3.11.} {\it If $\mathcal {K}_f$ is nonempty,
then we have that either $\mathcal{K}_f=\mathbb{N}^*$ or
$\mathcal{K}_f=\{1, ..., l\}$, where $l$ is an integer
such that $a^2(l+1)^2=D$.}
\begin{proof}
If $\mathcal{K}_f=\mathbb{N}^*$, then Lemma 3.11 is true.
If $\mathcal{K}_f\ne \mathbb{N}^*$, then the set
$\mathbb{N}^*\setminus\mathcal{K}_f$ is nonempty. By
the well-ordering principle (see, for example, page 13 of \cite{[A]}),
we know that $\mathbb{N}^*\setminus\mathcal{K}_f$
contains a smallest member, named
$s_0=\min (\mathbb{N}^*\setminus\mathcal{K}_f)$.
So $s_0\not\in{\mathcal K}_f$.
Suppose that $s_0=1$. Then $1\not\in{\mathcal K}_f$
and so $a^2=D$. It infers that ${\mathcal K}_f$ is empty, which
is impossible since $\mathcal{K}_f$ is nonempty.
Thus we have that $s_0\ge 2$ and all the
integers $s$ with $s<s_0$ are belonging to ${\mathcal K}_f$, i.e.,
$\{1,...,s_0-1\}\subseteq{\mathcal K}_f$.

On the other hand, since $s_0\not\in{\mathcal K}_f$, there is
an integer $s'$ with $1\le s'\le s_0$ such that $a^2s'^2-D=0$.
Thus $s'\in \mathbb{N}^*\setminus\mathcal{K}_f$. By
the minimality of $s_0$, one has $s'=s_0$. Hence $a^2s_0^2-D=0$,
which implies that $s_0+j\not\in{\mathcal K}_f$ for all nonnegative
integers $j$. Therefore ${\mathcal K}_f=\{1,...,s_0-1\}$ as desired.
Hence Lemma 3.11 is proved.
\end{proof}

Lemma 3.7 tells us that the sequence $\{d_{p^e}\}_{e=0}^\infty $
is nondecreasing. The following result describes
a condition on $f$ which guarantees that
$\{d_{p^e}\}_{e=0}^\infty $ is not a constant sequence.\\

\noindent {\bf Lemma 3.12.} {\it Let $\mathcal {K}_f$ be nonempty
and $k\in \mathcal {K}_f$. Then for any prime $p$,
there exists a unique nonnegative integer $e$ such that}
$d_{p^e}\le k<d_{p^{e+1}}.$
\begin{proof}
Let $\mathcal {K}_f$ be nonempty and $p$ be a prime. For any
$k\in \mathcal {K}_f$, we define the subset $R_p(k)$ of
$\mathbb{N}^*$ by $R_p(k):=\{i\in \mathbb{N}^*: k<d_{p^i}\}.$

If the set $R_p(k)$ is nonempty, then by the well-ordering principle
\cite{[A]}, we know that $R_p(k)$ contains
a smallest element, named $i_0=\min (R_p(k))$. Letting $e=i_0-1$ gives the
desired result $d_{p^e}\le k<d_{p^{e+1}}$. The uniqueness of $e$ follows
from Lemma 3.7. Thus we only need to show that $R_p(k)$ is nonempty
for any prime $p$. In the following we show the equivalent statement that
there is a positive integer $t$ such that $k<d_{p^{t}}$.

First let $p$ be a prime such that $p\mid a$. If $p|b$, then one can take $t=1$
since by Lemma 3.8, we have $d_{p^0}\le k< d_p=\infty$.
If $p\nmid b$, then pick $t$ to be a positive integer such that
$k<p^{t}$. However, Lemma 3.8 tells us that $d_{p^{t}}=p^{t}$.
Hence $k<d_{p^{t}}$. The statement is true for the case $p\mid a$.
Consequently, we let $p$ be a prime such that $p\nmid a$.

If $p=2$ and $D_4\not\equiv1\pmod8$, then we take
$t=2\lfloor\frac{\nu _2(D)}{2}\rfloor+1$. Then by Lemma 3.9 (ii),
we obtain that $d_{2^t}=\infty$ and thus $k<d_{2^t}$. The statement is true
in this case.

If $p$ is an odd prime with either $\nu _p(D)$ being odd or $(\frac{D_p}{p})=-1$, then
we take $t=\nu _p(D)+1$. But $d_{p^{\nu _p(D)+1}}=\infty$ by Lemma 3.10 (ii). So we get
the desired result $k<d_{p^t}$. The statement is proved in this case.

If either $p=2$ and $D_4\equiv1\pmod8$, or $p$ is an odd prime with $\nu _p(D)$
being even and $(\frac{D_p}{p})=1$, then by the assumption that
$\mathcal {K}_f$ is nonempty and Lemma 3.11, we have either
$\mathcal{K}_f=\mathbb{N}^*$, or $\mathcal{K}_f=\{1, ..., l\}$
for some positive integer $l$ satisfying $a^2(l+1)^2=D$. We take
\begin{align}\label{2.10}
t=\max\big(\nu _p(D)+1, \log_p(a^2k^2+|D|)+1\big)
\end{align}
if $\mathcal{K}_f=\mathbb{N}^*$, and
\begin{align}\label{2.11}
t=\max\big(\nu _p(D)+1, \max_{1\le i\le l}\{\nu _p(a^2i^2-D)\}\big)+1
\end{align}
if $\mathcal{K}_f$ is a finite set. So $t>\nu _p(D)$. It then follows from
Lemma 3.9 (iii) and Lemma 3.10 (iii) that $d_{p^{t}}$ is the smallest positive root of
$a^2x^2-D\equiv 0\pmod {p^{t+1}}$ if $p=2$ and $d_{p^{t}}$ is the smallest
positive root of $a^2x^2-D\equiv 0\pmod {p^{t}}$ if $p\ne 2$, respectively.
Now we show that $d_{p^t}>k$.

For the former case $\mathcal{K}_f=\mathbb{N}^*$, we have that $a^2i^2-D\ne 0$
for any $i\in \mathbb{N}^*$. By (\ref{2.10}), we have $t>\nu _p(D)$ and
$p^{t}>a^2k^2+|D|$. Since $a^2d_{p^{t}}^2-D\equiv 0\pmod {p^{t}}$, we can write
$a^2d_{p^{t}}^2-D=p^{t}u$ for some integer $u$. Then $u\ne 0$ because $a^2i^2-D\ne 0$
for any $i\in \mathbb{N}^*$. It then follows that $a^2d_{p^{t}}^2\ge p^{t}-|D|>a^2k^2$,
which implies that $d_{p^{t}}> k$ as required.

For the latter case that $\mathcal{K}_f=\{1, ..., l\}$, we have $a^2(l+1)^2-D=0$.
But by (\ref{2.11}), we have $a^2i^2-D\not\equiv 0\pmod {p^t}$ for all $1\le i\le l$.
Hence $l+1$ is the smallest positive root of the congruences
$a^2x^2-D\equiv0\pmod{p^{t}}$ and $a^2x^2-D\equiv0\pmod{p^{t+1}}$.
But $d_{p^t}$ is equal to the smallest positive root
of the congruence $a^2x^2-D\equiv0\pmod{p^{t+1}}$ if $p=2$,
and $d_{p^t}$ is equal to the smallest positive root
of the congruence $a^2x^2-D\equiv0\pmod{p^{t}}$ if $p\ne 2$.
Thus $d_{p^t}=l+1$. However, $k\le l$ since
$k\in {\mathcal K}_f$. So we obtain that $d_{p^{t}}>k$ as desired.
The proof of Lemma 3.12 is complete.
\end{proof}

\noindent {\bf Lemma 3.13.} {\it Let $\mathcal {K}_f$ be nonempty.

{\rm (i).} Let $a$ be odd and $D_4\equiv 1\pmod 8$. If $k\in \mathcal {K}_f$
and there is an integer $e>\nu _2(D)$ such that $d_{2^e}\le k<d_{2^{e+1}}$, then
$\max_{1\le i\le k }\{\nu _2(a^2i^2-D)\}=\nu _2(a^2d_{2^e}^2-D)=e+1.$

{\rm (ii).} Let $p$ be an odd prime with $p\nmid
a$. For any $k\in \mathcal {K}_f$ such that $d_{p^e}\le k<d_{p^{e+1}}$
for some nonnegative integer $e$, we have}
$\max_{1\le i\le k}\{\nu _p(a^2i^2-D)\}=\nu _p(a^2d_{p^e}^2-D)=e.$
\begin{proof}
(i).  By Lemma 3.9 (iii), $d_{2^{e}}$ is the smallest solution of
$a^2x^2-D\equiv0\pmod{2^{e+1}}$. Since $k\ge d_{2^e}$, we have
$\max_{1\le i\le k}\{\nu _2(a^2i^2-D)\}\ge \nu _2(a^2d_{2^e}^2-D)\ge e+1.$
Part (iii) of Lemma 3.9 also tells us that $d_{2^{e+1}}$ is the smallest
positive solution of $a^2x^2-D\equiv0\pmod{2^{e+2}}$. Thus $\nu _2(a^2i^2-D)<e+2$
for all $1\le i<d_{2^{e+1}}$. It then follows immediately from $k<d_{2^{e+1}}$ that
$\max_{1\le i\le k}\{\nu _2(a^2i^2-D)\}\le  e+1.$
Therefore we obtain the desired result
$\max_{1\le i\le k}\{\nu _2(a^2i^2-D)\}=\nu _2(a^2d_{2^e}^2-D)=e+1$.
Part (i) is proved.

(ii). Since $ k \ge d_{p^e}$, $d_{p^e}<\infty$.
By parts (i) and (iii) of Lemma 3.10, we derive that
$e\le \nu _p(a^2d_{p^e}^2-D)\le  \max_{1\le i\le k}\{\nu _p(a^2i^2-D)\}.$
Noticing the facts that $k<d_{p^{e+1}}$ and $d_{p^{e+1}}$ is
the smallest solution of $a^2x^2-D\equiv0\pmod{p^{e+1}}$, we obtain
$\max_{1\le i\le k}\{\nu _p(a^2i^2-D)\}\le e.$
It then follows that
$\max_{1\le i\le k}\{\nu _p(a^2i^2-D)\}= \nu _p(a^2d_{p^e}^2-D)=e$
as required. Part (ii) is true. So Lemma 3.13 is proved.
\end{proof}

\section {\bf $p$-Adic analysis of $g_{k,f}$ and determination of local periods}

In this section, we supply detailed local analysis to all eventually periodic
arithmetic functions $g_{k,f}$ presented in Section 2.
We always assume that $f$ is primitive throughout this section. Let
$S_{k,f}(n):=\{f(n), f(n+1), ..., f(n+k)\}$
for any positive integer $n$. By the definition of $d_{p^{e+1}}$, we know that
if $e$ is a nonnegative integer such that $d_{p^e}\le k<d_{p^{e+1}}$,
then there is at most one term divisible by $p^{e+1}$ in the set $S_{k,f}(n)$
for any positive integer $n$. We have
\begin{align}\label{4.1}
&g_{p,k,f}(n)= \sum_{m\in S_{k, f}(n)}v_{p}(m)-
{\rm max}_{m\in S_{k,f}(n)}\{v_{p}(m)\}\\
\nonumber&= \sum_{i=1}^\infty \#\{m\in S_{k,f}(n):\nu _p( m)\ge i\}-\sum_{i=1}^\infty (1 \
{\rm if} \ \nu _p(m)\ge i \
{\rm for \ some} \ m\in S_{k,f}(n))\\
\nonumber&= \sum_{i=1}^\infty \#\{m\in S_{k,f}(n):p^i\mid m\}-\sum_{i=1}^\infty (1 \
{\rm if} \ p^i \
{\rm divides \ some} \ m\in S_{k,f}(n))=\sum_{i=1}^\infty h_{p,i}(n),
\end{align}
where $h_{p,i}(n):=\max(0, \#\{m\in S_{k,f}(n): p^i\mid m\}-1)$.
It follows that if the set $\{m\in S_{k,f}(n): p^i\mid m\}$ is nonempty, then
$h_{p,i}(n):=\#\{m\in S_{k,f}(n): p^i\mid m\}-1.$

\noindent{\bf Lemma 4.1.} {\it  Let $\mathcal {K}_f$ be nonempty and $p$
be a prime with $p\nmid a$. Let $k\in \mathcal {K}_f$ and $e$ be the
positive integer such that $d_{p^{e}}\le k< d_{p^{e+1}}$. If $e\le \nu _p(D)$,
then there is a positive integer $n_0$ such that
$g_{p,k,f}(n_0+p^{\lceil e/2\rceil -1})\ne g_{p,k,f}(n_0)$.}
\begin{proof}
Let first $n$ be any positive integer and $e\le \nu _p(D)$. Consider the difference
$\Delta_1(n):=g_{p,k,f}(n+p^{\lceil e/2\rceil-1})-g_{p,k,f}(n).$
Then to show Lemma 4.1,
it suffices to find some suitable integer $n$ such that $\Delta_1(n)\ne 0$,
which will be done in the following.

For any integer $i\ge e+1$, since $d_{p^e}\le k<d_{p^{e+1}}$, there is at most one
term divisible by $p^{i}$ in the set $S_{k,f}(n)$ for any positive integer $n$.
Thus $\#\{m\in S_{k, f}(n): p^i\mid f(m)\}\le 1$ and so $h_{p, i}(n)=0$
for any positive integer $n$. It follows from (\ref{4.1}) that
\begin{align}\label{4.2}
\Delta_1(n)=\sum_{i=1}^{e}\big(h_{p, i}(n+p^{\lceil e/2\rceil-1})
-h_{p,i}(n)\big).
\end{align}
We claim that for any integers $m$ and $i$ with $1\le i\le e$, we have
\begin{align}\label{4.3}
\nu _p(f(n))\ge i \Longleftrightarrow \nu _p(f(n+mp^{\lceil e/2\rceil})\ge i.
\end{align}

For the case $p=2$ and $2\nmid a$, since $\nu _2(D)\ge e\ge 1$, we have that $\nu _2(D)\ge 2$
and $b$ is even.  If $\nu _2(f(n))\ge i$, then
$\Big(an+\frac{b}{2}\Big)^2\equiv \frac{D}{4}\pmod {2^i}.$
It follows that $ \nu _2(an+\frac{b}{2})\ge  \lceil
\frac{i}{2}\rceil $ if $i\le \nu _2(D)-2$, and $ \nu _2(an+\frac{b}{2})= \lceil
\frac{i}{2}\rceil-1$ if $i=\nu _2(D)-1\ \text{or}\ \nu _2(D)$. Hence for any
integer $m$, we obtain
\begin{align*}
\nu _2(f(n+m2^{\lceil
 e/2\rceil}))&=\nu _2(f(n)+(an+\frac{b}{2})m2^{\lceil
 e/2\rceil+1}+am^2\cdot2^{2\lceil  e/2\rceil})\\
&\ge \min\Big\{i,  \Big\lceil\frac{i}{2}\Big\rceil-1+\Big\lceil
\frac{e}{2}\Big\rceil+1+\nu _2(m), 2\Big\lceil \frac{e}{2}\Big\rceil +2\nu _2(m)\Big\}\ge i.
\end{align*}
Similarly, one has
\begin{align}\label{4.4}
\nu _2(f(n-m2^{\lceil e/2\rceil}))\ge i.
\end{align}

Conversely, if $\nu _2(f(n+m2^{ \lceil \frac{e}{2}\rceil})\ge i$,
then we obtain by replacing $n$ with $n+m2^{\lceil e/2\rceil}$ in (\ref{4.4}) that
$\nu _2(f(n))\ge i.$ Therefore for any integers $m$ and $i$ with $1\le i\le e$, we have
\begin{align}\label{4.5}
\nu _2(f(n))\ge i \Longleftrightarrow \nu _2(f(n+m2^{\lceil e/2\rceil})\ge i.
\end{align}
That is, the claim is true for the case $p=2$ and $2\nmid a$.

For the case $p\ne 2$ and $p\nmid a$, if $\nu _p(f(n))\ge i$, then it
follows from $e\le \nu _p(D)$ and $(2an+b)^2\equiv D\pmod{p^i}$ that $\nu _p(2an+b)\ge \lceil i/2\rceil$,
which implies that
\begin{align}\label{4.6}
\nu _p(f(n\pm mp^{\lceil e/2\rceil}))\ge \min\{i, \lceil i/2\rceil+\lceil e/2\rceil+\nu _p(m), 2\lceil
e/2\rceil+2\nu _p(m)\}\ge i.
\end{align}
If $\nu _p(f(n+mp^{\lceil e/2\rceil}))\ge i$, replacing $n$ with $n+mp^{\lceil e/2\rceil}$ in (\ref{4.6}),
then $\nu _p(f(n))=\nu _p(f(n+mp^{\lceil e/2\rceil}-
mp^{\lceil e/2\rceil}))\ge i$. Hence (\ref{4.3}) holds in this case. The claim is proved.

Replacing $e$ by $2\lceil e/2\rceil-2$,
then (\ref{4.3}) gives that for any given $1\le i\le 2\lceil e/2\rceil-2$
and any $0\le j\le k$,
$\nu _p(f(n+j))\ge i \Longleftrightarrow \nu _p(f(n+j+p^{\lceil e/2\rceil-1})\ge i.$
Thus the number of terms divisible by $p^i$ in $S_{k,f}(n)$ is equal to that in
$S_{k, f}(n+p^{\lceil e/2\rceil-1})$ for each $1\le i\le 2\lceil  e/2\rceil-2$.
It implies that
$h_{p, i}(n+p^{\lceil  e/2\rceil-1})=h_{p, i}(n)$
for each $1\le i\le 2\lceil  e/2\rceil-2$. Therefore by (\ref{4.2}),
we derive that
\begin{align}\label{4.7}
\Delta_1(n)=\sum_{i=2\lceil e/2\rceil-1}^{e}\big(h_{p, i}
(n+p^{\lceil e/2\rceil-1})-h_{p,i}(n)\big).
\end{align}

Since $k\ge d_{p^e}$ implying that $d_{p^e}<\infty$,
by the definition of $d_{p^e}$ we know that $S(f, p^e)$ is nonempty. Define
$x_0:=\mbox{the smallest positive solution of the congruence} \ f(x)\equiv0\pmod{p^e}.$
Then by Lemma 3.2 (i)-(ii) and Lemma 3.3 (i),
any term divisible by $p^e$ in the quadratic sequence $\{f(m)\}_{m=1}^{\infty}$
must be of the form $f(x_0+tp^{\lceil  e/2\rceil})$ with $t\in \mathbb{N}$.
But $\lceil  e/2\rceil =\lceil  (e-1)/2\rceil $ if $e$ is even. Thus by Lemma 3.2
(i)-(ii) and Lemma 3.3 (i),
the terms divisible by $p^{e-1}$ are exactly the terms divisible by $p^e$ in
the quadratic progression $\{f(m)\}_{m=1}^{\infty}$. Hence
the terms divisible by $p^{e-1}$ are exactly the terms divisible by $p^e$ in
the set $S_{k,f}(n)$ (resp. $S_{k,f}(n+p^{\lceil  e/2\rceil-1})$)
if $e$ is even. Namely,
$\{\bar m\in S_{k,f}(m): p^{e-1}\mid \bar m\}=\{\bar m\in S_{k,f}(m): p^{e}\mid \bar m\}$
for $m=n, n+p^{\lceil e/2\rceil-1}$.
So $h_{p,e-1}(m)=h_{p,e}(m)$ for $m=n, n+p^{\lceil e/2\rceil-1}$, which implies that
$h_{p, e-1}(n+p^{\lceil e/2\rceil-1})-h_{p,e-1}(n)=h_{p, e}(n+p^{\lceil e/2\rceil-1})-h_{p,e}(n)$
if $e$ is even. It then follows from (\ref{4.7}) that
\begin{align}\label{4.8}
\Delta_1(n)=2^{\frac{1+(-1)^e}{2}}\big(h_{p, e}
(n+p^{\lceil e/2\rceil-1})-h_{p,e}(n)\big).
\end{align}
Since the terms divisible by $p^e$ in the sets $S_{k,f}(n)$ and $S_{k,f}(n+p^{\lceil  e/2\rceil-1})$
are of the form $f(x_0+tp^{\lceil  e/2\rceil})$ with $t\in \mathbb{N}$,
in order to compute $\Delta_1(n)$, it is sufficient to compare
the number of terms of the form $f(x_0+tp^{\lceil e/2\rceil})$ ($t\in\mathbb{N}$)
in the set $S_{k,f}(n)$ with that in the set $S_{k,f}(n+p^{\lceil
 e/2\rceil-1})$. By Lemma 3.9 (i) and Lemma 3.10 (i), $d_{p^e}=p^{\lceil  e/2\rceil}$.
But $k\ge d_{p^e}$. Thus $k\ge p^{\lceil  e/2\rceil}$.
Since $\nu _p(k+1)<\lceil  e/2\rceil$, we can suppose that $k=k_0p^{\lceil e/2\rceil}+ r$
for unique two integers $k_0$ and $r$ with $k_0\ge 1$ and $0\leq r\leq p^{\lceil  e/2\rceil}-2$.

If $0\le r< p^{\lceil  e/2\rceil}-p^{\lceil  e/2\rceil-1}$, then $p^{\lceil e/2\rceil-1}\le r
+p^{\lceil e/2\rceil-1}<p^{\lceil e/2\rceil}$. Hence the number of integers $t$ such that
$x_0+p^{\lceil e/2\rceil-1}\le x_0+tp^{\lceil e/2\rceil}\le x_0+k+p^{\lceil e/2\rceil-1}$
is equal to $\Big\lfloor \frac{k+p^{\lceil e/2\rceil-1}}{p^{\lceil e/2\rceil}}\Big\rfloor
=k_0+\Big\lfloor \frac{r+p^{\lceil e/2\rceil-1}}{p^{\lceil e/2\rceil}}\Big\rfloor
=k_0.$
So there are exactly $k_0$ terms divisible by $p^e$ in the set $S_{k,f}(x_0+
p^{\lceil e/2\rceil-1})$. Thus $h_{p, e}(x_0+p^{\lceil e/2\rceil-1})=k_0-1$.
Similarly, by counting the number of integers
$t$ satisfying $x_0\le x_0+tp^{\lceil e/2\rceil}\le x_0+k$, we get that the number of terms
divisible by $p^e$ in $S_{k,f}(x_0)$ equals
$\lfloor \frac{k}{p^{\lceil e/2\rceil}}\rfloor+1=k_0+1$ and so $h_{p, e}(x_0)=k_0$.
Thus we derive from (\ref{4.8}) that $\Delta_1(x_0)=-2^{\frac{1+(-1)^e}{2}}.$

If $p^{\lceil  e/2\rceil}-p^{\lceil  e/2\rceil-1}\le r\le p^{\lceil  e/2\rceil}-2$, then
$p^{\lceil  e/2\rceil}-p^{\lceil  e/2\rceil-1}+1\le r+1\le p^{\lceil  e/2\rceil}-1$
and $p^{\lceil e/2\rceil}+1\le r+p^{\lceil  e/2\rceil-1}+1\le p^{\lceil
e/2\rceil}+p^{\lceil  e/2\rceil-1}-1.$
Therefore, by counting the number of integers $t$ such that
$x_0+1\le x_0+tp^{\lceil  e/2\rceil}\le x_0+k+1$ (resp.
$x_0+p^{\lceil  e/2\rceil-1}+1\le x_0+tp^{\lceil  e/2\rceil}
\le x_0+p^{\lceil  e/2\rceil-1}+k+1$), we deduce that
$h_{p, e}(x_0+1)=\Big\lfloor \frac{k+1}
{p^{\lceil e/2\rceil}}\Big\rfloor-1=k_0+\Big\lfloor \frac{r+1}
{p^{\lceil e/2\rceil}}\Big\rfloor-1=k_0-1$
and
$
h_{p,e}(x_0+p^{\lceil  e/2\rceil-1}+1)=\Big\lfloor \frac{k+p^{\lceil  e/2\rceil-1}+1}
{p^{\lceil e/2\rceil}}\Big\rfloor-1=k_0+\Big\lfloor \frac{r+p^{\lceil  e/2\rceil-1}+1}
{p^{\lceil e/2\rceil}}\Big\rfloor-1=k_0.
$
It then follows from (\ref{4.8}) that $\Delta_1(x_0+1)=2^{\frac{1+(-1)^e}{2}}.$
Thus the desired result follows immediately. This completes the proof of Lemma 4.1.
\end{proof}

With the help of (\ref{4.1}), we can make a detailed local analysis to
determine the local period $P_{p,k,f}$ for each prime factor $p$ of $B_k$. We
have the following results.

\noindent{\bf Lemma 4.2.} {\it Let $p$ be a prime such that $p|a$. Then}
\begin{align*}
P_{p,k,f}={\left\{
\begin{array}{rl}
p^{\nu _p(B_k)}, & {\it if}\  p\nmid b\  {\it and}\
\nu _p(k+1)<\nu _p(B_k),\\
1, & {\it otherwise}.\\
\end{array}
\right.}
\end{align*}
\begin{proof}
If $p|b$, then $p\nmid f(n)$ for any positive integer
$n$ since $\gcd(a,b,c)=1$. In other words, $g_{p,k,f}(n)=0$ for any
positive integer $n$. Thus $P_{p,k,f}=1$ as required. Lemma 4.2 is
true if $p|b$.

Now we let $p\nmid b$. Then $p\nmid D=b^2-4ac$ since $p|a$. It
follows that $\nu _p(a^2n^2-D)=0$ for any positive integer $n$. Hence
$\nu _p(B_k)=\nu _p({\rm lcm}_{1\le i\le k}\{i(a^2i^2-D)\})
=\max_{1\le i\le k}\{\nu _p(i)\}=\nu _p(L_k).$
By Lemma 3.1, there is exactly one term divisible by $p^{e}$ in any
consecutive $p^e$ terms of the quadratic progression $\{f(n+m)\}_{m\in \mathbb{N}}$
for any given positive integers $e$ and $n$. Since
$p^{\nu _p(L_k)}\le k< p^{\nu _p(L_k)+1}$, it follows from Lemma 3.8 that
$d_{p^{\nu _p(L_k)}}\le k<d_{p^{\nu _p(L_k)+1}}.$
Then there is at most one term divisible by $p^{\nu _p(L_k)+1}$ in $S_{k,f}(n)$
for any positive integer $n$. Consider the following two cases.

{\bf Case 1.} $\nu _p(k+1)\ge \nu _p(B_k)=\nu _p(L_k)$.  By Lemma 3.1, we deduce that there
are exactly $\frac{k+1}{p^e}$ terms divisible by $p^e$ in $S_{k,f}(n)$
(resp. $S_{k,f}(n+1)$) for any positive integer $n$ and each $e\in \{1,..., \nu _p(L_k)\}$.
On the other hand, since there is at most one term divisible by
$p^{\nu _p(L_k)+1}$ in $S_{k,f}(n)$ (resp. $S_{k,f}(n+1)$), we have by (\ref{4.1}) that
$g_{p,k,f}(n)=\sum_{e=1}^{\nu _p(L_k)}\Big(\frac{k+1}{p^e}-1\Big)=g_{p,k,f}(n+1)$
for any positive integer $n$. Therefore $P_{p,k,f}=1$ as desired.
Lemma 4.2 is proved in this case.

{\bf Case 2.} $\nu _p(k+1)< \nu _p(B_k)=\nu _p(L_k)$. Evidently, $\nu _p(L_k)\ge 1$.
Since there is at most
one term divisible by $p^{\nu _p(L_k)+1}$ in $S_{k,f}(n)$ for any positive
integer $n$, we have $h_{p, e}(n)=0$ if $e\ge \nu _p(L_k)+1$.
Thus we can deduce from (\ref{4.1}) that
$g_{p,k,f}(n)=\sum_{e=1}^{\nu _p(L_k)}h_{p,e}(n).$

By Lemma 2.2, $p^{\nu _p(L_k)}$ is a period of $g_{p,k,f}$.
So it remains to prove that $p^{\nu _p(L_k)-1}$ is not a period of
$g_{p,k,f}$. For any integer $e$ such that $1\le e\le \nu _p(L_k)-1$,
since $f(n+p^{\nu _p(L_k)-1})\equiv f(n)\pmod {p^e}$ for any positive
integer $n$, $p^{\nu _p(L_k)-1}$ is a period of $h_{p,e}$.
Hence we only need to prove that $p^{\nu _p(L_k)-1}$ is not a period of
$h_{p, \nu _p(L_k)}$. Since $\nu _p(k+1)< \nu _p(L_k)$, we can pick an $r\in
\{0,1,...,p^{\nu _p(L_k)}-2\}$ such that $k\equiv r\pmod {p^{\nu _p(L_k)}}$.

{\bf Subcase 2.1.} $0\le r< p^{\nu _p(L_k)}-p^{\nu _p(L_k)-1}$.
Then by Lemma 3.1, we can choose a positive
integer $n_0$ such that $f(n_0)\equiv 0 \pmod {p^{\nu _p(L_k)}}$.
And so the terms divisible by $p^{\nu _p(L_k)}$ in the quadratic sequence $\{
f(n_0+i)\}_{i\in \mathbb{N}}$ must be of the form
$f(n_0+tp^{\nu _p(L_k)})$ for some $t\in \mathbb{N}$. It then follows that there
are exactly $1+\Big\lfloor \frac{k}{p^{\nu _p(L_k)}}\Big\rfloor$ terms divisible by
$p^{\nu _p(L_k)}$ in $S_{k,f}(n_0)$ and there are exactly
$\Big\lfloor \frac{k+p^{\nu _p(L_k)-1}}{p^{\nu _p(L_k)}}\Big\rfloor=\Big\lfloor
\frac{k-r}{p^{\nu _p(L_k)}}\Big\rfloor+\Big\lfloor\frac{p^{\nu _p(L_k)-1}
+r}{p^{\nu _p(L_k)}}\Big\rfloor=\Big\lfloor \frac{k}{p^{\nu _p(L_k)}}\Big\rfloor$
terms divisible by $p^{\nu _p(L_k)}$ in
$S_{k,f}(n_0+p^{\nu _p(L_k)-1})$, where the last equality is derived from
$k\equiv r\pmod {p^{\nu _p(L_k)}}$ and $0\le r< p^{\nu _p(L_k)}-p^{\nu _p(L_k)-1}$. Thus
$h_{p, \nu _p(L_k)}(n_0)=\Big\lfloor\frac{k}{p^{\nu _p(L_k)}}\Big\rfloor=h_{p,
\nu _p(L_k)}(n_0+p^{\nu _p(L_k)-1})+1.$
That is, $p^{\nu _p(L_k)-1}$ is not a period of $h_{p, \nu _p(L_k)}$.

{\bf Subcase 2.2.} $p^{\nu _p(L_k)}-p^{\nu _p(L_k)-1}\le r\le p^{\nu _p(L_k)}-2$.
Again by Lemma 3.1, we can pick a suitable positive integer $m_0$ such that
$f(m_0+p^{\nu _p(L_k)-1}-1)\equiv 0\pmod {p^{\nu _p(L_k)}}.$
It follows that the terms divisible by $p^{\nu _p(L_k)}$ in the quadratic sequence $\{
f(m_0+i)\}_{i\in \mathbb{N}}$ must be of the form
$f(m_0+p^{\nu _p(L_k)-1}-1+sp^{\nu _p(L_k)})$ for some $s\in \mathbb{N}$.
Since $k\equiv r\pmod {p^{\nu _p(L_k)}}$ and $p^{\nu _p(L_k)}-p^{\nu _p(L_k)-1}\le r\le p^{\nu _p(L_k)}-2$,
we can derive that the number of terms divisible by $p^{\nu _p(L_k)}$ in
the set $S_{k,f}(m_0)$ is equal to
$1+\Big\lfloor\frac{k-(p^{\nu _p(L_k)-1}-1)}{p^{\nu _p(L_k)}}\Big\rfloor=1+\Big\lfloor
\frac{k}{p^{\nu _p(L_k)}}\Big\rfloor$
and the number of terms divisible by $p^{\nu _p(L_k)}$ in the set
$S_{k,f}(m_0+p^{\nu _p(L_k)-1})$ equals
$\Big\lfloor \frac{k+1}{p^{\nu _p(L_k)}}\Big\rfloor=\Big\lfloor
\frac{k}{p^{\nu _p(L_k)}}\Big\rfloor.$
Thus
$h_{p, \nu _p(L_k)}(m_0)=h_{p, \nu _p(L_k)}(m_0+p^{\nu _p(L_k)-1})+1.$
Namely, $p^{\nu _p(L_k)-1}$ is not a period of $h_{p, \nu _p(L_k)}$
as required. The proof of Lemma 4.2 is complete.
\end{proof}

Now we need only to handle the even prime $2$ and the
odd prime $p$ with $p\nmid a$, respectively. We first consider the case $2\nmid a$.
Since $D_4\equiv1\pmod4$ if $\nu _2(D)=0$, we have that $\nu _2(D)\ge 1$ if $e=\nu _2(D)$
with $D_4\not\equiv 1\pmod4$. Therefore, if either $e=\nu _2(D)$ with $D_4\not\equiv
1\pmod4$ or $e>\nu _2(D)$ with $D_4\not\equiv1\pmod8$, then by Lemma 3.9 (ii), $d_{2^e}=\infty$.
But there is no integer $k$ such that $k\ge d_{2^e}$ for such integers $e$.
So one only needs to consider the cases occurred exactly in Lemma 4.3.

\noindent{\bf Lemma 4.3.} {\it  Let  $a$ be odd and $\mathcal {K}_f$ be
nonempty. Let $k\in \mathcal {K}_f$ and $e$ be the nonnegative integer
such that $d_{2^{e}}\le k< d_{2^{e+1}}$. Each of the following is true.

{\rm (i).} If $e=\nu _2(D)$ with $D_4\equiv 1\pmod 4$ or $e\le  2\lfloor
\frac{\nu _2(D)}{2}\rfloor-1$, then
\begin{align*}
 P_{2,k,f}={\left\{
  \begin{array}{rl}
  2^{\lceil  e/2 \rceil}, & \mbox{if} \ \nu _2(k+1)<\lceil  e/2 \rceil,\\
1, & \mbox{if} \ \nu _2(k+1)\ge \lceil  e/2 \rceil.
 \end{array}
\right.}
\end{align*}

{\rm (ii).} If $e>\nu _2(D)$ with $D_4\equiv1\pmod8$, then
$P_{2,k,f}=2^{e-\frac{\nu _2(D)}{2}}$. }
\begin{proof}
Since $d_{2^e}\le k<d_{2^{e+1}}$, there is at
most one term divisible by $2^{e+1}$ in $S_{k,f}(n)$ for any positive
integer $n$. It follows from (\ref{4.1}) that $2^0=1$ is the
smallest period of $g_{2,k,f}$ if $e=0$. So it remains to treat with the case
$e\ge 1$. Let now $e\ge 1$ and $n\ge 1$ be an arbitrary given integer.
Since $\#\{m\in S_{k,f}(n): 2^i\mid
m\}\le 1$ if $i\ge e+1$, then by (\ref{4.1}),
\begin{align}\label{4.11}
g_{2,k,f}(n)=\sum_{i=1}^{e}h_{2, i}(n),
\end{align}
where
$$h_{2, i}(n)={\rm max} (0, \#\{m\in S_{k,f}(n): 2^i\mid m\}-1)=
{\rm max} (0, \#\{0\le j\le k: 2^i\mid f(n+j)\}-1).$$
Clearly, $h_{2, i}(n)=\#\{0\le j\le k: 2^i\mid f(n+j)\}-1$ if there
is at least one term divisible by $2^i$ in $S_{k,f}(n)$.

(i).  Since $e=\nu _2(D)\ge 1$ with $D_4\equiv5\pmod8$
or $e\le 2\lfloor \frac{\nu _2(D)}{2}\rfloor-1$, we have $\nu _2(D)\ge 2$
and $b$ is even. If $ \nu _2(k+1)\ge \lceil  e/2\rceil$, comparing $S_{k,f}(n)$ with
$S_{k,f}(n+1)$, we find that their distinct terms are $f(n)$ and
$f(n+k+1)$. Since $ \nu _2(k+1)\ge \lceil  e/2\rceil$, we have
$k+1=m_02^{\lceil  e/2\rceil}$ for some positive integer $m_0$.
From (\ref{4.5}), we deduce that for any given integer $i$ with $1\le i\le e$,
$\nu _2(f(n))\ge i$ if and only if $\nu _2(f(n+k+1))\ge i$. Thus
the number of terms divisible by $2^i$ in $S_{k,f}(n)$ is
equal to the number of terms divisible by $2^i$ in $S_{k,f}(n+1)$
for each $i\in \{1,...,e\}$. Hence by (\ref{4.11}), we obtain that
$g_{2,k,f}(n)=g_{2,k,f}(n+1)$ for any positive integer $n$, which
implies that $P_{2,k,f}=1$. Part (i) is true in this case.

In what follows we let $\nu _2(k+1)<\lceil  e/2 \rceil$. It follows from (\ref{4.5})
that for any given $1\le i\le e$ and for any $0\le j\le k$,
$\nu _2(f(n+j))\ge i \Longleftrightarrow \nu _2(f(n+j+2^{\lceil e/2\rceil})\ge i.$
In other words, the number of terms divisible by $2^i$ in $S_{k,f}(n+2^{\lceil
 e/2\rceil})$ is equal to the number of terms divisible by
$2^i$ in $S_{k,f}(n)$ for each $1\le i\le e$. So
$g_{2,k,f}(n+2^{\lceil  e/2\rceil})=g_{2,k,f}(n)$ for any
positive integer $n$. This infers that $2^{\lceil  e/2\rceil}$
is a period of $g_{2,k,f}$. On the other hand, by Lemma 4.1
one knows that there is a positive integer $n_0$
such that $g_{2,k,f}(n_0+2^{\lceil  e/2\rceil-1})=g_{2,k,f}(n_0)$.
Thus $2^{\lceil  e/2\rceil-1}$ is not a period of $g_{2,k,f}$.
Therefore $2^{\lceil  e/2\rceil}$ is the smallest period
of $g_{2,k,f}$. Part (i) is proved.

(ii). Since $D_4$ is odd, $\nu _2(D)$ is even.
First, we prove that $2^{e-\frac{\nu _2(D)}{2}}$ is a period of $g_{2,k,f}$.
Since $f(m+2^e)\equiv f(m)\pmod {2^i}$, we get that $h_{2, i}(m+2^e)
=h_{2, i}(m)$ for any integers $m$ and $i$ with $0\le i\le e$. So by (\ref{4.2}),
$g_{2, k, f}(m+2^e)=g_{2, k, f}(m)$ for any integer $m$, i.e.,
$2^e=2^{e-\frac{\nu _2(D)}{2}}$ is a period of $g_{2, k, f}$ if
$\nu _2(D)=0$.

Now let $\nu _2(D)\ge 2$. Then $b$ is even.
Let $l$ be any given positive integer with $l\ge \nu _2(D)$,
and let $i\in \{1,...,l\}$ and $j\in\{0,...,k\}$. If
\begin{align}\label{4.12}
\nu _2(f(n+j))\ge i,
\end{align}
then $(a(n+j)+\frac{b}{2})^2\equiv \frac{D}{4}\pmod {2^i}$,
which implies that
$\nu _2(a(n+j)+\frac{b}{2})\ge \min\{\frac{\nu _2(D)}{2}-1, \lceil
\frac{i}{2}\rceil\}.$
But $\min\{\frac{\nu _2(D)}{2}-1, \lceil
\frac{i}{2}\rceil\}\ge\lceil \frac{i}{2}\rceil$ if
$i\le \nu _2(D)-2$, and $\min\{\frac{\nu _2(D)}{2}-1, \lceil
\frac{i}{2}\rceil\}=\frac{\nu _2(D)}{2}-1$ if
$i\ge \nu _2(D)-1$. It then follows that
\begin{align*}
 &\nu _2(f(n+j+2^{l-\frac{\nu _2(D)}{2}}))=\nu _2\Big(f(n+j)+\big(a(n+j)+\frac{b}{2}\big)
2^{l-\frac{\nu _2(D)}{2}+1}+a2^{2l-\nu _2(D)}\Big)\\
&\ge \min\Big\{i, l-\frac{\nu _2(D)}{2}+1+\min\big\{\frac{\nu _2(D)}{2}-1, \big\lceil
\frac{i}{2}\big\rceil\big\}, 2l-\nu _2(D)\Big\} \ge i.
\end{align*}
Similarly, we have
\begin{align}\label{4.13}
\nu _2(f(n+j-2^{l-\frac{\nu _2(D)}{2}}))\ge i.
\end{align}

If $\nu _2(f(n+j+2^{l-\frac{\nu _2(D)}{2}}))\ge i$, then the process of (\ref{4.13}) derived
from (\ref{4.12}) with $n$ replaced by $n+2^{l-\frac{\nu _2(D)}{2}}$ gives us that $\nu _2(f(n+j))\ge i$.
Therefore, if $l$ is an integer with $l\ge \nu _2(D)$, then for any integers $i$ and $j$ with
$1\le i\le l$ and $0\le j\le k$, we have
\begin{align}\label{4.14}
\nu _2(f(n+j))\ge i\Longleftrightarrow \nu _2(f(n+j+2^{l-\frac{\nu _2(D)}{2}}))\ge i.
\end{align}
Since $e> \nu _2(D)$, the number of terms divisible by $2^{i}$ in
$S_{k,f}(n+2^{e-\frac{\nu _2(D)}{2}})$ is equal to the number of terms
divisible by $2^{i}$ in $S_{k,f}(n)$ and
so $h_{2, i}(n)=h_{2, i}(n+2^{e-\frac{\nu _2(D)}{2}})$ for each $1\le i\le e$.
Thus by (\ref{4.11}), we have $g_{2,k,f}(n+2^{e-\frac{\nu _2(D)}{2}})=g_{2,k,f}(n)$
for any positive integer $n$. So $2^{e-\frac{\nu _2(D)}{2}}$ is a period
of $g_{2,k,f}$.

In the following, we prove that $2^{e-\frac{\nu _2(D)}{2}-1}$ is not the period of
$g_{2,k,f}$. It suffices to show that $\Delta_2(n)\ne 0$ for some integer $n$, where
\begin{align}\label{4.15}
\Delta_2(n):=g_{2,k,f}(n+2^{e-\frac{\nu _2(D)}{2}-1})-g_{2,k,f}(n)
=\sum_{i=1}^e\big(h_{2,i}(n+2^{e-\frac{\nu _2(D)}{2}-1})-h_{2,i}(n)\big).
\end{align}
Since $e-1\ge \nu _2(D)$, replacing $e$ by $e-1$ in (\ref{4.14}), one gets that for any integers $i$ and $j$ with
$1\le i\le e-1$ and $0\le j\le k$, $\nu _2(f(n+j))\ge i$ if and only if
$\nu _2(f(n+j+2^{e-\frac{\nu _2(D)}{2}-1}))\ge i$.  In other words, the number of terms
divisible by $2^i$ in $S_{k,f}(n)$ is equal to that in
$S_{k,f}(n+2^{e-\frac{\nu _2(D)}{2}-1})$ for each $1\le i\le e-1$. Thus
$\sum_{i=1}^{e-1}h_{2,i}(n)=\sum_{i=1}^{e-1}h_{2,i}
(n+2^{e-\frac{\nu _2(D)}{2}-1}).$
It then follows from (\ref{4.15}) that $\Delta_2(n)=h_{2,e}(n+2^{e-\frac{\nu _2(D)}{2}-1})-h_{2,e}(n).$
Therefore, our final task is to find some suitable integer $n$ such that
\begin{align}\label{4.16}
h_{2,e}(n+2^{e-\frac{\nu _2(D)}{2}-1})\ne h_{2,e}(n).
\end{align}

Since $a$ is odd and $D_4\equiv1\pmod8$, we have
$\nu _2(a^2(2^{\frac{\nu _2(D)}{2}})^2-D)=\nu _2(D)+\nu _2(a^2-D_4)\ge
\nu _2(D)+3,$
which means that $a^2(2^{\frac{\nu _2(D)}{2}})^2-D\equiv0\pmod{2^{\nu _2(D)+3}}.$
We can easily check that $2^{\frac{\nu _2(D)}{2}}$ is the smallest solution of the
congruence $a^2x^2-D\equiv0\pmod{2^i}$ for any $i$ with $\nu _2(D)\le i\le \nu _2(D)+3$.
Then by parts (i) and (iii) of Lemma 3.8, we derive that
$d_{2^{\nu _2(D)-1}}=d_{2^{\nu _2(D)}}=d_{2^{\nu _2(D)+1}}=d_{2^{\nu _2(D)+2}}=2^{\frac{\nu _2(D)}{2}}.$
Since $d_{2^e}\le k<d_{2^{e+1}}$ and $e> \nu _2(D)$, we have $e\ge \nu _2(D)+2$.
But $D_4\equiv1\pmod 8$, then by part (iv) of Lemma 3.2, $S(f, 2^e)$ is nonempty.

We claim that for any $i\ge e$, $d_{2^i}<2^{i-\frac{\nu _2(D)}{2}-1}.$
In fact, by Lemma 3.9 (iii), $d_{2^i}$ equals the smallest positive
solution of $a^2x^2-D\equiv0\pmod{2^{i+1}}$. Then $\nu _2(d_{2^i})=\frac{\nu _2(D)}{2}$.
So $\nu _2(a^2(2^{i-\frac{\nu _2(D)}{2}}-d_{2^i})^2-D)\ge i+1$, and hence $2^{i-\frac{\nu _2(D)}{2}}-d_{2^i}$
is also a solution of $a^2x^2-D\equiv0\pmod{2^{i+1}}$. From the minimality of $d_{2^i}$, we get that
$2^{i-\frac{\nu _2(D)}{2}}-d_{2^i}\ge d_{2^i}$. But $2^{i-\frac{\nu _2(D)}{2}}-d_{2^i}\ne d_{2^i}$. Otherwise,
$\nu _2(d_{2^i})=i-\frac{\nu _2(D)}{2}-1\ge \nu _2(D)+2-\frac{\nu _2(D)}{2}-1=\frac{\nu _2(D)}{2}+1,$
which is a contradiction since $\nu _2(d_{2^i})=\frac{\nu _2(D)}{2}$. So
$d_{2^i}<2^{i-\frac{\nu _2(D)}{2}-1}$. The claim is proved. From the claim, we know that
$d_{2^e}<2^{e-\frac{\nu _2(D)}{2}-1}$ and $d_{2^{e+1}}<2^{e-\frac{\nu _2(D)}{2}}.$
Thus $k<d_{2^{e+1}}<2^{e-\frac{\nu _2(D)}{2}}$.

If either $d_{2^{e+1}}> 2^{e-\frac{\nu _2(D)}{2}-1}$ and $k<2^{e-\frac{\nu _2(D)}{2}-1}$
or $d_{2^{e+1}}\le 2^{e-\frac{\nu _2(D)}{2}-1}$, then $d_{2^e}\le k< 2^{e-\frac{\nu _2(D)}{2}-1}$.
Since $S(f, 2^e)$ is nonempty, by Lemma 3.6 we can choose a positive integer $n_0$ such that
$\nu _2(f(n_0))\ge e$ and $\nu _2(f(n_0+d_{2^e}))\ge e.$ By Lemma 3.2 (iii), the
terms divisible by $2^e$ in the quadratic progression
$\{f(n_0+i)\}_{i\in \mathbb{N}}$ must be of the form $f(n_0+t_1
2^{e-\frac{\nu _2(D)}{2}})$ or $f(n_0+d_{2^e}+t_2
2^{e-\frac{\nu _2(D)}{2}})$, where $t_1,t_2\in \mathbb{N}$. On the one hand, since
$d_{2^e}\le k< 2^{e-\frac{\nu _2(D)}{2}-1}$,
$f(n_0)$ and $f(n_0+d_{2^e})$ are the only two terms  divisible by $2^e$
in $S_{k,f}(n_0)$. On the other hand, since
$n_0+d_{2^e}<n_0+2^{e-\frac{\nu _2(D)}{2}-1}+j\le n_0+
k+2^{e-\frac{\nu _2(D)}{2}-1}<n_0+2^{e-\frac{\nu _2(D)}{2}}$
for all $0\le j\le k$, there is no term divisible by $2^e$ in the set
$S_{k,f}(n_0+2^{e-\frac{\nu _2(D)}{2}-1})$. It follows that
$h_{2, e}(n_0+2^{e-\frac{\nu _2(D)}{2}-1})=0$
and $h_{2, e}(n_0)=1$. So (\ref{4.16}) is true in this case.

If $k\ge 2^{e-\frac{\nu _2(D)}{2}-1}$, then it follows from Lemma 3.6 and the fact that
$S(f, 2^e)$ is nonempty that there is a positive integer $n_1$ so that
$\nu _2(f(n_1+2^{e-\frac{\nu _2(D)}{2}-1}-1-d_{2^e}))\ge e$
and $\nu _2(f(n_1+2^{e-\frac{\nu _2(D)}{2}-1}-1))\ge e.$
Hence Lemma 3.2 (iii) tells us that the terms divisible by $2^e$
in the quadratic progression $\{f(n_1+i)\}_{i\in \mathbb{N}}$ should
be of the form
\begin{align}\label{4.17}
f(n_1+2^{e-\frac{\nu _2(D)}{2}-1}-1-d_{2^e}+t_1 2^{e-\frac{\nu _2(D)}{2}})
\end{align}
or
\begin{align}\label{4.18}
f(n_1+2^{e-\frac{\nu _2(D)}{2}-1}-1+t_2 2^{e-\frac{\nu _2(D)}{2}}),
\end{align}
where $t_1,t_2\in \mathbb{N}$. Since $k<d_{2^{e+1}}$ and $d_{2^{e+1}}<2^{e-\frac{\nu _2(D)}{2}}$,
we have $2^{e-\frac{\nu _2(D)}{2}-1} \le
j+2^{e-\frac{\nu _2(D)}{2}-1}<2^{e-\frac{\nu _2(D)}{2}}+2^{e-\frac{\nu _2(D)}{2}-1}-1$
for all $j\in \{0,1,..., k\}$. Therefore, there is at most one term of the form (\ref{4.17})
with $t_1\in \mathbb{N}$ and no term of the form (\ref{4.18}) with $t_2\in \mathbb{N}$
in the set $S_{k,f}(n_1+2^{e-\frac{\nu _2(D)}{2}-1})$.
Since $k\ge 2^{e-\frac{\nu _2(D)}{2}-1}$ and $d_{2^e}<2^{e-\frac{\nu _2(D)}{2}-1}$ ,
$f(n_1+2^{e-\frac{\nu _2(D)}{2}-1}-1)$ and
$f(n_1+2^{e-\frac{\nu _2(D)}{2}-1}-1-d_{2^e})$ are the only two terms
divisible by $2^e$ in the set $S_{k,f}(n_1)$. So $h_{2,e}(n_1)=1$ and
$h_{2,e}(n_1+2^{e-\frac{\nu _2(D)}{2}})=0$, which implies that (\ref{4.16}) holds in the case.
This concludes that $2^{e-\frac{\nu _2(D)}{2}}$ is the smallest period of $g_{2,k,f}$.

This completes the proof of Lemma 4.3.
\end{proof}

In what follows, we treat with all the odd primes $p$ with $p\nmid 2a$ and $p| B_k$.

\noindent{\bf Lemma 4.4.} {\it Let $\mathcal {K}_f$ be nonempty and
$p$ be an odd prime with $p\nmid a$. Let $k\in \mathcal {K}_f$ and $e$
be the nonnegative integer such that $d_{p^e}\le k<d_{p^{e+1}}$. Then}

\begin{align*}
 P_{p,k,f}={\left\{
  \begin{array}{rl}
  1, & {\it if \ either} \  e\le \nu _p(D)\ {\it and}\ \nu _p(k+1)\ge \lceil e/2\rceil,\\
  & {\it or} \  e> \nu _p(D)\ {\it and}\ \nu _p(k+1)\ge e-\nu _p(D)/2,\\
p^{\lceil e/2\rceil},&  {\it if} \ e\le \nu _p(D) \ {\it and} \ \nu _p(k+1)<\lceil e/2\rceil,\\
p^{e-\nu _p(D)/2}, & {\it if} \ e>\nu _p(D)\  {\it and}\ \nu _p(k+1)<e-\nu _p(D)/2.
 \end{array}
\right.}
\end{align*}
\begin{proof}
Let $n\ge 1$ be any positive integer. Since $d_{p^e}\le k<d_{p^{e+1}}$, there is at most one
term divisible by $p^{e+1}$ in the set $S_{k,f}(n)$ for any
positive integer $n$. It follows from (\ref{4.1}) that
\begin{align}\label{4.19}
g_{p,k,f}(n)=\sum_{i=1}^e h_{p, i}(n),
\end{align}
where $h_{p, i}(n)=\#\{0\le j\le k: p^i\mid f(n+j)\}-1$
if there is at least one term divisible by $p^i$ in $S_{k,f}(n)$.
Otherwise, $h_{p, i}(n)=0$. Thus $g_{p,k,f}(n)=0$ for any
positive integer $n$, and so $P_{p,k,f}=1$ if $e=0$.

In what follows we let $e\ge 1$. Note that if $e>\nu _p(D)$ and $d_{p^e}<\infty$, then by
parts (ii) and (iii) of Lemma 3.10, we know that $\nu _p(D)$ is even
and $(D_p/p)=1$ for such primes $p$.

First we show that if $l\ge \nu _p(D)$ is an integer, then for any
integers $m$ and $i$ with $1\le i\le l$, we have
\begin{align}\label{4.20}
\nu _p(f(n))\ge i \Longleftrightarrow \nu _p(f(n+mp^{l-\frac{\nu _p(D)}{2}})\ge i.
\end{align}
In fact, by (\ref{4.3}), we know that for any integers $m$ and $i$ with $1\le i\le \nu _p(D)$,
\begin{align}\label{4.21}
\nu _p(f(n))\ge i \Longleftrightarrow \nu _p(f(n+mp^{\nu _p(D)})\ge i.
\end{align}
Since $l\ge \nu _p(D)$, $p^{l-\frac{\nu _p(D)}{2}}$ is a multiple of
$p^{\frac{\nu _p(D)}{2}}$. Then by (\ref{4.21}), (\ref{4.20})
is true for each $1\le i\le \nu _p(D)$. For $\nu _p(D)\le i\le l$, we can deduce from
$f(n)\equiv0\pmod {p^i}$ that $\nu _p(2an+b)\ge \nu _p(D)/2$, which implies that
\begin{align}\label{4.22}
\nu _p(f(n\pm mp^{l-\frac{\nu _p(D)}{2}}))
\ge \min\{i, \nu _p(D)/2+l-\nu _p(D)/2+\nu _p(m), 2l-\nu _p(D)+\nu _p(m)\}\ge i.
\end{align}
Conversely, if $\nu _p(f(n+mp^{l-\frac{\nu _p(D)}{2}}))\ge i$, then replacing $n$
with $n+mp^{l-\frac{\nu _p(D)}{2}}$ in (\ref{4.22}), we get
$\nu _p(f(n))=\nu _p(f(n+mp^{l-\frac{\nu _p(D)}{2}}-mp^{l-\frac{\nu _p(D)}{2}}))\ge i$.
Hence (\ref{4.20}) is proved.

If either $e\le \nu _p(D)$ and $ \nu _p(k+1)\ge \lceil e/2\rceil$, or $e> \nu _p(D)$ and
$\nu _p(k+1)\ge e-\nu _p(D)/2$, then either $p^{\lceil e/2\rceil}\mid (k+1)$,
or $p^{e-\nu _p(D)/2}\mid (k+1)$. It then follows immediately from (\ref{4.3}) and (\ref{4.20})
with $l=e$ that for each $1\le i\le e$, $\nu _p(f(n))\ge i$ if and only if $\nu _p(f(n+k+1))\ge i$.
But the distinct terms of the sets $S_{k,f}(n)$ and $S_{k,f}(n+1)$ are $f(n)$ and
$f(n+k+1)$. Thus the number of terms divisible by $p^i$ in $S_{k,f}(n)$
is equal to that in $S_{k,f}(n+1)$ for any  $i\in \{1,..., e\}$.
Thus we have $h_{p, i}(n+1)=h_{p, i}(n)$ for each $i\in
\{1,...,e\}$, and so by (\ref{4.19}), $g_{p,k,f}(n+1)=g_{p,k,f}(n)$ for any
positive integer $n$. Hence $P_{p,k,f}=1$. So Lemma 4.4 is true if
either $e\le \nu _p(D)$ and $ \nu _p(k+1)\ge \lceil e/2\rceil$, or $e> \nu _p(D)$ and
$\nu _p(k+1)\ge e-\nu _p(D)/2$.

Now let $e\le \nu _p(D)$ and $\nu _p(k+1)<\lceil e/2\rceil$. Taking $m=1$ in (\ref{4.3}),
we have that for any given $1\le i\le e$ and for any $0\le j\le k$,
$\nu _p(f(n+j))\ge i$ if and only if $\nu _p(f(n+j+p^{\lceil e/2\rceil})\ge i$.
In other words, the number of terms divisible by $p^i$ in $S_{k,f}(n)$
is equal to that in $S_{k,f}(n+p^{\lceil e/2\rceil})$ for any  $i\in \{1,..., e\}$.
It infers that $h_{p, i}(n+p^{\lceil e/2\rceil})=h_{p, i}(n)$ for each
$i\in\{1,...,e\}$. Thus by (\ref{4.19}) $g_{p,k,f}(n+p^{\lceil
e/2\rceil})=g_{p,k,f}(n)$ for any positive integer $n$, and so
$p^{\lceil e/2\rceil}$ is a period of $g_{p,k,f}$.
But Lemma 4.1 implies that there is a positive integer $n_0$ such that
$g_{p,k,f}(n_0+p^{\lceil e/2\rceil-1})\ne g_{p,k,f}(n_0).$ Therefore
$p^{\lceil e/2\rceil-1}$ is not the period of $g_{p,k,f}$.
Thus $p^{\lceil e/2\rceil}$ is the smallest period of $g_{p,k,f}$ as required.
Thus Lemma 4.4 is true if $e\le \nu _p(D)$ and $\nu _p(k+1)<\lceil e/2\rceil$.

We only need to deal with the remaining case: $e> \nu _p(D)$ and $\nu _p(k+1)<e-\nu _p(D)/2$
which will be done in what follows. First, from (\ref{4.20}) with $l=e$ and $m=1$,
it follows that for any given $1\le i\le e$ and for any
$0\le j\le k$, $\nu _p(f(n+j))\ge i$ if and only if $\nu _p(f(n+j+p^{e-\nu _p(D)/2}))\ge i$.
Namely, the number of terms divisible by $p^i$ in $S_{k,f}(n)$
is equal to that in $S_{k,f}(n+p^{e-\nu _p(D)/2})$ for each $i\in \{1,..., e\}$.
Hence $g_{p,k,f}(n+p^{e-\nu _p(D)/2})=g_{p,k,f}(n)$ for any positive integer
$n$ by (\ref{4.19}). Thus $p^{e-\nu _p(D)/2}$ is a period of $g_{p,k,f}$.

By Lemma 3.3 (iii), we know that the congruence
$f(x)\equiv0\pmod{p^e}$ has exactly two solutions in the interval
$[1, p^{e-\frac{\nu _p(D)}{2}}]$. It follows that $d_{p^e}\le
(p^{e-\nu _p(D)/2}-1)/2$. Therefore, we can find a positive
integer $u_0$ with $1\le u_0\le \frac{p+1}{2}$ such that
$(u_0-1)p^{e-\nu _p(D)/2-1}\le d_{p^e}<u_0p^{e-\nu _p(D)/2-1}.$
To prove that $p^{e-\nu _p(D)/2}$ is the smallest period of
$g_{p,k,f}$, it suffices to
prove that $u_0p^{e-\frac{\nu _p(D)}{2}-1}$ is not a period of $g_{p,k,f}$.
For this purpose, we define the arithmetic function $\Delta $ for any
positive integer $n$ by
\begin{align}\label{4.23}
\Delta(n):=g_{p,k,f}(n+u_0p^{e-\nu _p(D)/2-1})-g_{p,k,f}(n).
\end{align}

Since $e>\nu _p(D)$, we have $e-1\ge \nu _p(D)$. Picking $l=e-1$ and $m=u_0$ in (\ref{4.20}),
we get that for any given $1\le i\le e-1$ and for any
$0\le j\le k$, $\nu _p(f(n+j+u_0p^{e-\nu _p(D)/2-1}))\ge i$
if and only if $\nu _p(f(n+j))\ge i$. Hence the number of terms divisible
by $p^i$ in $S_{k,f}(n)$ is equal to that in $S_{k,f}(n+u_0p^{e-\nu _p(D)/2-1})$,
i.e., $h_{p,i}(n+u_0p^{e-\nu _p(D)/2-1})=h_{p,i}(n)$ for each $1\le i\le e-1$.
So by (\ref{4.19}) and (\ref{4.23}), we get
\begin{align}\label{4.24}
\Delta(n)=h_{p,e}(n+u_0p^{e-\nu _p(D)/2-1})-h_{p,e}(n).
\end{align}

Define the two sets
$\mathcal{A}_1(n):=\{f(n),...,f(n+d_{p^e}),...,f(n+u_0p^{e-\nu _p(D)/2-1}-1)\}
$ and $\mathcal{A}_2(n):=\{f(n+k+1),...,f(n+k+u_0p^{e-\nu _p(D)/2-1})\}.$
Evidently, $S_{k,f}(n+u_0p^{e-\nu _p(D)/2-1})\subseteq \mathcal {A}_2(n)$
if $k<u_0p^{e-\nu _p(D)/2-1}$. If $k\ge u_0p^{e-\nu _p(D)/2-1}$,
then we have the following disjoint unions:
$S_{k,f}(n)=\mathcal{A}_1(n) \bigcup  \{f(n+u_0p^{e-\nu _p(D)/2-1}),...,f(n+k)\}
$ and $S_{k,f}(n+u_0p^{e-\nu _p(D)/2-1})= \{f(n+u_0p^{e-\nu _p(D)/2-1}),...,f(n+k)\}
 \bigcup \mathcal{A}_2(n).$
Claim that there is a positive integer $n_0$ such that the set $S_{k,f}(n_0)$
contains exactly two terms divisible by $p^e$ if $k<u_0p^{e-\nu _p(D)/2-1}$,
while the set $\mathcal{A}_1(n_0)$ holds exactly two terms divisible by $p^e$
and the set $\mathcal{A}_2(n_0)$ has at most one term divisible by $p^e$.

Suppose that the claim is true. If $k<u_0p^{e-\nu _p(D)/2-1}$, then it follows from the claim that
$
h_{p, e}(n_0+u_0p^{e-\nu _p(D)/2-1})=\max(0, \#\{ m\in S_{k,f}(n_0+u_0p^{e-\nu _p(D)/2-1}): p^e\mid m\}-1 )=0
$
and
$
h_{p, e}(n_0)= \#\{ m\in S_{k,f}(n_0): p^e\mid m\}-1=1.
$
Hence by (\ref{4.24}), we get $\Delta(n_0)=-1$.
If $k\ge u_0p^{e-\nu _p(D)/2-1}$, then we derive from the claim that
$
h_{p, e}(n_0+u_0p^{e-\nu _p(D)/2-1})\le \#\{u_0p^{e-\nu _p(D)/2-1}\le j\le k: p^e\mid f(n_0+j) \}
$
and
$
h_{p, e}(n_0)= \#\{u_0p^{e-\nu _p(D)/2-1}\le j\le k: p^e\mid f(n_0+j)\}+1.
$
It follows from (\ref{4.24}) that $\Delta(n_0)\le -1$. Therefore
$u_0p^{e-\nu _p(D)/2-1}$ is not a period of $g_{p,k,f}$.
Thus Lemma 4.4 is true if $e> \nu _p(D)$ and $\nu _p(k+1)<e-\nu _p(D)/2$.
It remains to prove that the claim is true.

First note that by Lemma 3.3 (iii), there are exactly two terms divisible by
$p^e$ in any consecutive $p^{e-\nu _p(D)/2}$ terms of the quadratic progression
$\{f(n)\}_{n=1}^{\infty}$. Since $\nu _p(k+1)<e-\nu _p(D)/2$, we can find some integer $r$
with $1\le r\le p^{e-\nu _p(D)/2}-1$ such that $k+1\equiv r\pmod{p^{e-\nu _p(D)/2}}.$
We divide the proof of the claim into the following two cases.

{\bf Case 1.} $r\in [1, u_0p^{e-\frac{\nu _p(D)}{2}-1}]$ with
$u_0\in[1, \frac{p-1}{2}]$, or $r\in [1, (p-1)p^{e-\nu _p(D)/2-1}/2]\cup (d_{p^e},
(p+1)p^{e-\nu _p(D)/2-1}/2]$ with $u_0=\frac{p+1}{2}$.
By Lemma 3.6 we can choose a positive integer
$n_0$ such that $\nu _p(f(n_0))\ge e$ and $\nu _p(f(n_0+d_{p^e}))\ge e$.
By Lemma 3.3 (iii), we know that the terms divisible
by $p^e$ in the quadratic progression $\{f(n_0+j)\}_{j\in
\mathbb{N}}$ must be of the form $f(n_0+t_1p^{e-\nu _p(D)/2})$ or
$f(n_0+d_{p^e}+t_2p^{e-\nu _p(D)/2})$, $t_1,t_2\in \mathbb{N}$.
Since $|\mathcal{A}_1(n_0)|=u_0p^{e-\nu _p(D)-1}<p^{e-\nu _p(D)/2}$
and $d_{p^e}<u_0p^{e-\nu _p(D)/2-1}$, $f(n_0)$ and $f(n_0+d_{p^e})$ are the
exactly two terms divisible by $p^e$ in $\mathcal{A}_1(n_0)$.
On the other hand, since $k\ge d_{p^e}$ and
$|S_{k,f}(n_0)|=k+1\le u_0p^{e-\nu _p(D)/2-1}<p^{e-\nu _p(D)/2},$
$f(n_0)$ and $f(n_0+d_{p^e})$ are exactly the two terms divisible by $p^e$ in $S_{k,f}(n_0)$
if $k<u_0p^{e-\nu _p(D)/2-1}$. Namely, $\mathcal{A}_1(n_0)$ holds exactly two terms divisible by $p^e$
and $S_{k,f}(n_0)$ contains exactly two terms divisible by $p^e$ if $k<u_0p^{e-\nu _p(D)/2-1}$.
Now we show that $\mathcal{A}_2(n_0)$ has at most
one term divisible by $p^e$. Since $|\mathcal{A}_2(n_0)|=u_0p^{e-\nu _p(D)/2-1}<p^{e-\nu _p(D)/2}$, there is
at most one term of the form $f(n_0+t_1p^{e-\nu _p(D)/2})$ and there is
at most one term of the form  $f(n_0+d_{p^e}+t_2p^{e-\nu _p(D)/2})$ in the set
$\mathcal{A}_2(n_0)$ with $t_1,t_2\in \mathbb{N}$. Therefore, we only need to show that
either there is no term of the form $f(n_0+t_1p^{e-\nu _p(D)/2})$, or there is
no term of the form  $f(n_0+d_{p^e}+t_2p^{e-\nu _p(D)/2})$ in the set
$\mathcal{A}_2(n_0)$, where $t_1,t_2\in \mathbb{N}$,
which will be done in the following.

If $r\in [1, u_0p^{e-\frac{\nu _p(D)}{2}-1}]$ with
$u_0\in [1, \frac{p-1}{2}]$, we have that for all $1\leq j\le u_0p^{e-\nu _p(D)/2-1}$,
$
k+j\equiv r+j-1\not \equiv 0\pmod {p^{e-\nu _p(D)/2}}
$
since
$1\le r+j-1\le 2u_0p^{e-\nu _p(D)/2-1}-1 < p^{e-\nu _p(D)/2}-1.$
Hence there is no term of  the
form $f(n_0+t_1p^{e-\nu _p(D)/2})$ in the set
$\mathcal{A}_2(n_0)$ with $t_1\in
\mathbb{N}$.

If $r\in [1, (p-1)p^{e-\nu _p(D)/2-1}/2]\cup
(d_{p^e}, (p+1)p^{e-\nu _p(D)/2-1}/2]$ with $u_0=\frac{p+1}{2}$, we have for
all $1\leq j\le u_0p^{e-\nu _p(D)/2-1}$ that
$1\le r+j-1\le (p-1)p^{e-\nu _p(D)/2-1}/2+u_0p^{e-\nu _p(D)/2-1}-1=p^{e-\nu _p(D)/2}-1$
if $r\in [1, (p-1)p^{e-\nu _p(D)/2-1}/2]$ and that
$d_{p^e}< r+j-1\le (p+1)p^{e-\nu _p(D)/2-1}-1\le p^{e-\nu _p(D)/2}+d_{p^e}-1$
if $r\in (d_{p^e}, (p+1)p^{e-\nu _p(D)/2-1}/2]$ since
$d_{p^e}\ge (u_0-1)p^{e-\nu _p(D)/2-1}=\frac{p-1}{2}p^{e-\nu _p(D)/2-1}\ge p^{e-\nu _p(D)/2-1}.$
That is, for all $1\leq j\le u_0p^{e-\nu _p(D)/2-1}$, we have
$
k+j\equiv r+j-1\not \equiv 0\pmod {p^{e-\nu _p(D)/2}}
$
if $r\in [1, (p-1)p^{e-\nu _p(D)/2-1}/2]$ and
$
k+j\equiv r+j-1\not \equiv d_{p^e}\pmod {p^{e-\nu _p(D)/2}}
$
if $r\in (d_{p^e}, (p+1)p^{e-\nu _p(D)/2-1}/2]$.
Therefore, there is no term of the form $f(n_0+t_1p^{e-\nu _p(D)})$ in
$\mathcal{A}_2(n_0)$ if $r\in [1, (p-1)p^{e-\nu _p(D)/2-1}/2]$ with $u_0=\frac{p+1}{2}$,
and there is no term of the form $f(n_0+d_{p^e}+t_2p^{e-\nu _p(D)/2-1})$
in $\mathcal{A}_2(n_0)$ if $r\in (d_{p^e}, (p+1)p^{e-\nu _p(D)/2-1}/2]$ with $u_0=\frac{p+1}{2}$,
where $t_1, t_2\in \mathbb{N}$. So the claim is proved for Case 1.

{\bf Case 2.} $r\in (u_0p^{e-\frac{\nu _p(D)}{2}-1}, p^{e-\nu _p(D)/2}-1]$ with
$u_0\in[1, \frac{p-1}{2}]$, or
$r\in ((p+1)p^{e-\frac{\nu _p(D)}{2}-1}/2,\\ p^{e-\nu _p(D)/2}-1]\cup((p-1)p^{e-\nu _p(D)/2-1}/2,d_{p^e}]$
with $u_0=\frac{p+1}{2}$. Then by Lemma 3.6, we can select a positive
integer $n_0$ such that $\nu _p(f(n_0+u_0p^{e-\nu _p(D)/2-1}-1))=e$ and
$\nu _p(f(n_0+u_0p^{e-\nu _p(D)/2-1}-1-d_{p^e}))\ge e$. Now by Lemma 3.3 (iii),
the terms divisible by $p^e$ in the quadratic progression
$\{f(n_0+j)\}_{j\in \mathbb{N}}$ are of the form
\begin{align}\label{4.28}
f(n_0+u_0p^{e-\nu _p(D)/2-1}-1+t_1p^{e-\nu _p(D)/2})
\end{align}
or
\begin{align}\label{4.29}
f(n_0+u_0p^{e-\nu _p(D)/2-1}-1-d_{p^e}+t_2p^{e-\nu _p(D)/2}),
\end{align}
where $t_1,t_2\in \mathbb{N}$.

Since $k+1\equiv r\pmod{p^{e-\nu _p(D)/2}}$, one may let $k=r-1+tp^{e-\nu _p(D)}$
for some integer $t\ge 0$. It follows that $k\ge r-1\ge u_0p^{e-\nu _p(D)/2-1}$
if $r\in (u_0p^{e-\frac{\nu _p(D)}{2}-1}, p^{e-\nu _p(D)/2}-1]$ with
$u_0\in[1, \frac{p-1}{2}]$, or
$r\in ((p+1)p^{e-\frac{\nu _p(D)}{2}-1}/2, p^{e-\nu _p(D)/2}-1]$
with $u_0=\frac{p+1}{2}$. If $r\in ((p-1)p^{e-\nu _p(D)/2-1}/2,d_{p^e}]$
with $u_0=\frac{p+1}{2}$, then it follows from $k\ge d_{p^e}$ and $r\le d_{p^e}$ that
$t\ge 1$. Thus $k\ge r-1+p^{e-\nu _p(D)}\ge u_0p^{e-\nu _p(D)/2-1}$. That is, we always have
$k\ge u_0p^{e-\nu _p(D)/2-1}$ in Case 2. Hence to finish the proof of the claim for
Case 2, we only need to treat with the two sets $\mathcal{A}_1(n_0)$ and $\mathcal{A}_2(n_0)$.

Evidently, $u_0p^{e-\nu _p(D)/2-1}-1-d_{p^e}\ge 0$ since
$d_{p^e}<u_0p^{e-\nu _p(D)/2-1}$. Again using the fact $|\mathcal{A}_1(n_0)|
=u_0p^{e-\nu _p(D)-1}<p^{e-\nu _p(D)/2}$, we know that $f(n_0+u_0p^{e-\nu _p(D)/2-1}-1)$ and
$f(n_0+u_0p^{e-\nu _p(D)/2-1}-1-d_{p^e})$ are the exactly two terms divisible by
$p^e$ in $\mathcal{A}_1(n_0)$.

Since $|\mathcal{A}_2(n_0)|=u_0p^{e-\nu _p(D)/2-1}$, to show that $\mathcal{A}_2(n_0)$ has at most
one term divisible by $p^e$, it is enough to show that either there is
no term of the form (\ref{4.28}), or there is no term of the form  (\ref{4.29}) in the set
$\mathcal{A}_2(n_0)$, where $t_1,t_2\in \mathbb{N}$.

If either $r\in (u_0p^{e-\frac{\nu _p(D)} {2}-1}, p^{e-\nu _p(D)/2}-1]$ with $u_0\in [1, (p-1)/2]$,
or $r\in ((p+1)p^{e-\frac{\nu _p(D)}{2}-1}/2, p^{e-\nu _p(D)/2}-1]$ with $u_0=\frac{p+1}{2}$,
then for all $1\le j\le u_0p^{e-\nu _p(D)/2-1}$, we have
$u_0p^{e-\nu _p(D)/2-1}<r+j-1 < p^{e-\nu _p(D)/2}+ u_0p^{e-\nu _p(D)/2-1}-1,$
which implies that
$k+j\equiv r+j-1\not\equiv u_0p^{e-\nu _p(D)/2-1}-1\pmod{p^{e-\nu _p(D)/2}}.$
Hence there is no term of the form (\ref{4.28}) with $t_1\in\mathbb{N}$ in $\mathcal{A}_2(n_0)$.

If $r\in ( (p-1)p^{e-\nu _p(D)/2-1}/2,d_{p^e}]$ with $u_0=\frac{p+1}{2}$, then
\begin{align*}
\frac{p-1}{2}p^{e-\nu _p(D)/2-1}&<r+j-1 \le  d_{p^e}+\frac{p+1}{2}p^{e-\nu _p(D)/2-1}-1\\
&<p^{e-\nu _p(D)/2}+\frac{p+1}{2}p^{e-\nu _p(D)/2-1}-1-d_{p^e}
\end{align*}
for all $1\le j\le u_0p^{e-\nu _p(D)/2-1}$ since $d_{p^e}\le \frac{p^{e-\nu _p(D)/2}-1}{2}$.
However,
\begin{align*}
\frac{p+1}{2}p^{e-\nu _p(D)/2-1}-1-d_{p^e}
&<\frac{p+1}{2}p^{e-\nu _p(D)/2-1}-1-\frac{p-1}{2}p^{e-\nu _p(D)/2-1}\\
&=p^{e-\nu _p(D)/2-1}-1<\frac{p-1}{2}p^{e-\nu _p(D)/2-1}.
\end{align*}
It then follows that
$k+j\equiv r+j-1\not\equiv \frac{p+1}{2}p^{e-\nu _p(D)/2-1}-1-d_{p^e}\pmod{p^{e-\nu _p(D)/2}}$
for all $1\le j\le u_0p^{e-\nu _p(D)/2-1}$. So there is no term of the form (\ref{4.29})
with $t_2\in\mathbb{N}$ in the set $\mathcal{A}_2(n_0)$.
So the claim is true in Case 2.

The proof of Lemma 4.4 is complete.
\end{proof}

From Lemmas 4.3 and 4.4, we see that $\nu _p(P_{p,k,f})$
depends on some nonnegative integer $e$ satisfying
$d_{p^e}\le k<d_{p^{e+1}}$. In other words,  we still don't get the
explicit value of $P_{p,k,f}$. Thus, to determine the exact value of
$P_{p,k,f}$ for those primes $p\nmid a$, we need to transform the
information on $e$ into explicit information depending on $k$ and $f$.
We have the following results.

\noindent{\bf Lemma 4.5.} {\it Let $a$ be odd and $\mathcal {K}_f$ be nonempty.
Then for any $k\in \mathcal {K}_f$, we have}
\begin{align*}
P_{2,k,f}={\left\{
  \begin{array}{rl}
  2^{\nu _2(B_k)-2\nu _2(L_k)}, & {\it if} \  k<2^{\lfloor \frac{\nu _2(D)}
  {2}\rfloor} \ {\it and} \ \nu _2(k+1)<\nu _2(L_k),\\
  2^{\lfloor \frac{\nu _2(D)}{2}\rfloor}, & {\it if} \  k\ge
  2^{\lfloor \frac{\nu _2(D)}{2}\rfloor},\ D_4\not\equiv1\pmod8
  \ {\it and} \ \nu _2(k+1)<\lfloor \frac{\nu _2(D)}{2}\rfloor,\\
2^{\nu _2(B_k)-\nu _2(D)-1},&  {\it if} \  k\ge
  2^{\lfloor \frac{\nu _2(D)}{2}\rfloor} \ {\it and} \
  D_4\equiv1 \pmod 8,\\
1, & {\it otherwise.}
 \end{array}
\right.}
\end{align*}

\begin{proof}
Since $\mathcal {K}_f$ is nonempty and $k\in \mathcal {K}_f$,
by Lemma 3.12, there is a unique nonnegative integer $e$
such that $d_{2^e}\le k<d_{2^{e+1}}$. Consider the following three cases.

{\bf Case 1.} $k<2^{\lfloor\frac{\nu _2(D)}{2}\rfloor}$. Since $2\nmid a$, we have
$\nu _2(a^2i^2)=2\nu _2(i)\le 2(\lfloor \frac{\nu _2(D)}{2}\rfloor-1)\le \nu _2(D)-2<\nu _2(D)$
and so $\nu _2(a^2i^2-D)=2\nu _2(i)$ for any integer $i$ with $1\le i\le k$. Hence
\begin{align}\label{4.30}
\max_{1\le i\le k}\{ \nu _2(a^2i^2-D)\}=\max_{1\le i\le k}\{2\nu _2(i)\}=2\nu _2(L_k).
\end{align}
But

\begin{align}\label{4.31}
\nu _2(B_k)=\max_{1\le i\le k}\{\nu _2(i(a^2i^2-D))\}
=\max_{1\le i\le k}\{\nu _2(i)+\nu _2(a^2i^2-D)\} =\max_{1\le i\le
k}\{3\nu _2(i)\}=3\nu _2(L_k).
\end{align}
By Lemma 3.9 (i), we have $d_{2^{2\lfloor\frac{\nu _2(D)}{2}\rfloor-1}}
=2^{\lfloor \frac{\nu _2(D)}{2}\rfloor}$.
Thus by Lemma 3.7 and $d_{2^e}\le k<2^{\lfloor \frac{\nu _2(D)}{2}\rfloor}$, we have
$e<2\lfloor\frac{\nu _2(D)}{2}\rfloor-1.$
Notice that by part (i) of Lemma 3.9, $d_{2^e}$ (resp. $d_{2^{e+1}}$)
is the smallest positive root of the congruence
$a^2x^2-D\equiv 0\pmod{2^e}$ (resp. $a^2x^2-D\equiv 0\pmod{2^{e+1}}$).
Hence $2^{e+1}\nmid (a^2l^2-D)$ for all positive integers
$l<d_{2^{e+1}}$. But $d_{2^e}\le k< d_{2^{e+1}}$. Thus
$\max_{1\le i\le k}\{ \nu _2(a^2i^2-D)\}\ge e\ \mbox{and}\ \max_{1\le i\le k}\{ \nu _2(a^2i^2-D)\}<e+1.$
Then by (\ref{4.30}), $e=\max_{1\le i\le k}\{ \nu _2(a^2i^2-D)\}=2\nu _2(L_k).$
Therefore, by Lemma 4.3 (i) and (\ref{4.31}), we get that
$P_{2,k,f}=2^{\lceil  e/2\rceil}=2^{\nu _2(L_k)}=2^{\nu _2(B_k)-2\nu _2(L_k)}$
if $k<2^{\lfloor \frac{\nu _2(D)}{2}\rfloor}$ and $\nu _2(k+1)<\nu _2(L_k)$,
and $P_{2,k,f}=1$ if $k<2^{\lfloor \frac{\nu _2(D)}{2}\rfloor}$ and
$\nu _2(k+1)\ge \nu _2(L_k)$. Thus Lemma 4.5 is true in this case.

{\bf Case 2.} $k\ge 2^{\lfloor \frac{\nu _2(D)}{2}\rfloor}$ and
$D_4\not\equiv1\pmod8$. By parts (i) and (ii) of Lemma 3.9, one knows that
$d_{2^{2\lfloor \frac{\nu _2(D)}{2}\rfloor-1}}=2^{\lfloor\frac{\nu _2(D)}{2}\rfloor}$
and $d_{2^{2\lfloor \frac{\nu _2(D)}{2}\rfloor}}=\infty$ if $D_4\not\equiv1\pmod4$,
and $d_{2^{\nu _2(D)}}=2^{\frac{\nu _2(D)}{2}}$ and $d_{2^{\nu _2(D)+1}}=\infty$
if $D_4\equiv5\pmod8$.
It then follows from $d_{2^e}\le k<d_{2^{e+1}}$ and
$k\ge 2^{\lfloor \frac{\nu _2(D)}{2}\rfloor}$ that
$$
e={\left\{
  \begin{array}{rl}
 2\lfloor \frac{\nu _2(D)}{2}\rfloor-1, & \mbox{if} \ D_4\not\equiv1\pmod4,\\
 \nu _2(D), & \mbox{if}\ D_4\equiv5\pmod8.
 \end{array}
\right.}
$$
Thus by Lemma 4.3 (i), we obtain that
$P_{2,k,f}=2^{\lceil e/2\rceil}=2^{\lfloor\frac{\nu _2(D)}{2}\rfloor}$
if $k\ge 2^{\lfloor \frac{\nu _2(D)}{2}\rfloor}$, $D_4\not\equiv1\pmod8$
and $\nu _2(k+1)<\lfloor \frac{\nu _2(D)}{2}\rfloor$,
and $P_{2,k,f}=1$ if $k\ge 2^{\lfloor \frac{\nu _2(D)}{2}\rfloor}$,
$D_4\not\equiv1\pmod8$ and $\nu _2(k+1)\ge \lfloor \frac{\nu _2(D)}{2}\rfloor$.
Lemma 4.5 is true in Case 2.

{\bf Case 3.} $k\ge 2^{\lfloor \frac{\nu _2(D)}{2}\rfloor}$ and $D_4\equiv1\pmod8$.
Then $\nu _2(D)$ is even and $\nu _2(L_k)\ge \frac{\nu _2(D)}{2}$. Since $2^{\frac{\nu _2(D)}{2}}$
is equal to the smallest positive root of the congruence $a^2x^2-D\equiv0\pmod {2^{\nu _2(D)+3}}$,
we derive  from Lemma 3.9 (iii) that $d_{2^{\nu _2(D)+2}}=2^{\frac{\nu _2(D)}{2}}$.
Hence from Lemma 3.7 and $d_{2^e}\le k<d_{2^{e+1}}$, we can derive that $e\ge \nu _2(D)+2$.
So Lemma 3.13 (i) gives that $e=\max_{1\le i\le k}\{\nu _2(a^2i^2-D)\}-1.$
It follows from Lemma 4.3 (ii) that $\nu _2(P_{2,k,f})=e-\frac{\nu _2(D)}{2}
=\max_{1\le i\le k}\{\nu _2(a^2i^2-D)\}-\frac{\nu _2(D)}{2}-1.$
Therefore, to show that $\nu _2(P_{2,k,f})=\nu _2(B_k)-\nu _2(D)-1$, it suffices to
prove that the following is true:
\begin{align}\label{4.32}
\nu _2(B_k)=\max_{1\le i\le k}\{\nu _2(a^2i^2-D)\}+\frac{\nu _2(D)}{2},
\end{align}
which will be done in what follows.
Let $i$ be an integer such that $1\le i\le k$. Then we have that
\begin{align}\label{4.33}
\nu _2(i(a^2i^2-D))=\nu _2(i)+\min(2\nu _2(i), \nu _2(D))<\nu _2(L_k)+\nu _2(D)
\end{align}
if $\nu _2(i)<\frac{\nu _2(D)}{2}$, and that
\begin{align}\label{4.34}
\nu _2(i(a^2i^2-D))=\nu _2(i)+\nu _2(a^2i^2-D)=\frac{\nu _2(D)}{2}+\nu _2(a^2i^2-D)
\end{align}
if $\nu _2(i)=\frac{\nu _2(D)}{2}$, and that
\begin{align}\label{4.35}
\nu _2(i(a^2i^2-D))=\nu _2(i)+\min\{ \nu _2(a^2i^2), \nu _2(D)\}
=\nu _2(i)+\nu _2(D)\le \nu _2(L_k)+\nu _2(D)
\end{align}
if $\nu _2(i)>\frac{\nu _2(D)}{2}$. Now we claim that
\begin{align}\label{4.36}
\max_{1\le i\le k\atop \nu _2(i)=\nu _2(D)/2}\{\nu _2(a^2i^2-D)\}\ge \nu _2(L_k)+\frac{\nu _2(D)}{2}+1.
\end{align}
This is equivalent to showing that there is an integer
$i_0\in [1, k]$ with $\nu _2(i_0)=\frac{\nu _2(D)}{2}$ such that
$\nu _2(a^2i_0^2-D)\}\ge \nu _2(L_k)+\frac{\nu _2(D)}{2}+1.$

If $\nu _2(L_k)\le \frac{\nu _2(D)}{2}+2$ and $D_4\equiv1\pmod8$,
then pick $i_0=2^{\frac{\nu _2(D)}{2}}\in [1, k]$. Since $2\nmid a$, we have
$\nu _2(a^2(2^{\frac{\nu _2(D)}{2}})^2-D)= \nu _2(D)+\nu _2(a^2-D_4)\ge \nu _2(D)+3\ge \nu _2(L_k)+\frac{\nu _2(D)}{2}+1.$

If $\nu _2(L_k)>\frac{\nu _2(D)}{2}+2$ and $D_4\equiv1\pmod8$, then $\nu _2(L_k)+\frac{\nu _2(D)}{2}+1>\nu _2(D)+3$.
Since the discriminant of $a^2x^2-D$ is $4a^2D$ and $\nu _2(4a^2D)=\nu _2(D)+2$, then Lemma 3.2 (iv)
applied to the congruence $a^2x^2-D\equiv 0\pmod {2^{\nu _2(L_k)+\frac{\nu _2(D)}{2}+1}}$, we can find
an integer $i_0\in [1, 2^{\nu _2(L_k)}]\subseteq [1, k]$ satisfying that
$\nu _2(i_0)=\frac{\nu _2(D)}{2}$ and $a^2i_0^2-D\equiv0 \pmod{2^{\nu _2(L_k)+\frac{\nu _2(D)}{2}+1}}.$
The claim (\ref{4.36}) is proved. It follows from (\ref{4.33})-(\ref{4.36}) that
$\nu _2(B_k)=\frac{\nu _2(D)}{2}+\max_{1\le i\le k\atop \nu _2(i)=\nu _2(D)/2}\{\nu _2(a^2i^2-D)\}.$
On the other hand, since $\nu _2(a^2i^2-D)<\nu _2(D)$ if $\nu _2(i)<\frac{\nu _2(D)}{2}$ and
$\nu _2(a^2i^2-D)=\nu _2(D)$ if $\nu _2(i)>\frac{\nu _2(D)}{2}$, we have
$\max_{1\le i\le k}\{\nu _2(a^2i^2-D)\}
=\max_{1\le i\le k\atop \nu _2(i)=\nu _2(D)/2}\{\nu _2(a^2i^2-D)\}.$
Hence (\ref{4.32}) follows immediately. Lemma 4.5 is true for Case 3.

This ends the proof of Lemma 4.5.
\end{proof}

\noindent {\bf Lemma 4.6.} {\it Let $\mathcal {K}_f$ be nonempty. Then
for any $k\in \mathcal {K}_f$ and any odd prime $p$ with $p\nmid a$, we have}
\begin{align*}
P_{p,k,f}={\left\{
  \begin{array}{rl}
  p^{\nu _p(B_k)-2\nu _p(L_k)}, & {\it if} \  k<p^{\lceil \frac{\nu _p(D)}
  {2}\rceil} \ {\it and }\ \nu _p(k+1)<\nu _p(L_k),\\
  p^{\lceil \frac{\nu _p(D)}  {2}\rceil}, & {\it if} \
  k\ge p^{\lceil \frac{\nu _p(D)}{2}\rceil},\
  \nu _p(k+1)<\lceil \frac{\nu _p(D)}{2}\rceil \\ &
  {\it and \ either} \ 2\nmid \nu _p(D) \  {\it or} \ (\frac{D_p}{p})=-1,\\
p^{\nu _p(B_k)-\nu _p(D)},&  {\it if} \  k\ge p^{\lceil \frac{\nu _p(D)}{2}\rceil},
  \nu _p(k+1)<\nu _p(B_k)-\nu _p(D), \\& 2|\nu _p(D)\  \ {\it and} \ (\frac{D_p}{p})=1,\\
1, & {\it otherwise}.
 \end{array}
\right.}
\end{align*}
\begin{proof}
By Lemma 3.12, we can find a unique nonnegative integer $e$ such that $d_{p^e}\le
k<d_{p^{e+1}}$ since $\mathcal{K}_f$ is nonempty and $k\in \mathcal{K}_f$.
Let $p$ be an odd prime with $p\nmid a$. Then
\begin{align}\label{4.39}
\nu _p(B_k)=\nu _p({\rm lcm}_{1\le i\le k}\{i(a^2i^2-D)\})=\max_{1\le
i\le k}\{\nu _p(i(a^2i^2-D))\}.
\end{align}

If $k<p^{\lceil \frac{\nu _p(D)}{2}\rceil}$, then for any integer
$i$ with $1\le i\le k$, we have
\begin{align}\label{4.40}
2\nu _p(i)\le 2\nu _p(L_k)\le 2(\lceil \frac{\nu _p(D)}{2}\rceil-1)\le \nu _p(D)-1,
\end{align}
which implies that $\nu _p(a^2i^2-D)=2\nu _p(i)$. Hence
\begin{align}\label{4.41}
\max_{1\le i\le k}\{\nu _p(a^2i^2-D)\}=\max_{1\le i\le k}\{2\nu _p(i)\}=2\nu _p(L_k)
\end{align}
and by (\ref{4.39}), we have $\nu _p(B_k)=\max_{1\le i\le k}\{\nu _p(i)+2\nu _p(i)\}=3\nu _p(L_k).$
Since $p\nmid a$, by Lemma 3.13 (ii), we have $e=\max_{1\le i\le k}\{\nu _p(a^2i^2-D)\}$.
It then follows from (\ref{4.40}) and (\ref{4.41}) that $e=2\nu _p(L_k)<\nu _p(D)$. Thus by Lemma 4.4,
$P_{p,k,f}=p^{\lceil  e/2\rceil}=p^{\nu _p(L_k)}=p^{\nu _p(B_k)-2\nu _p(L_k)}$
if $k<p^{\lceil\frac{\nu _p(D)}{2}\rceil}$ and $\nu _p(k+1)<\nu _p(L_k)$,
and $P_{p,k,f}=1$ if $k<p^{\lceil\frac{\nu _p(D)}{2}\rceil}$ and $\nu _p(k+1)\ge \nu _p(L_k)$.
So Lemma 4.6 is true if $k<p^{\lceil \frac{\nu _p(D)}{2}\rceil}$.
In what follows we let $k\ge p^{\lceil \frac{\nu _p(D)}{2}\rceil}$. Then $\nu _p(L_k)\ge
\lceil \frac{\nu _p(D)}{2}\rceil$.

If $2\nmid \nu _p(D)$ or $(\frac{D_p}{p})=-1$, then by parts (i) and (ii)
of Lemma 3.10, we have $d_{p^{\nu _p(D)}}=p^{\lceil
\frac{\nu _p(D)}{2}\rceil}$ and $d_{p^{\nu _p(D)+1}}=\infty$.
Since $k\ge p^{\lceil \frac{\nu _p(D)}{2}\rceil}$ and $d_{p^e}\le k<d_{p^{e+1}}$,
we obtain by Lemma 3.7 that $e=\nu _p(D)$. It follows from Lemma 4.4
that $P_{p,k,f}=p^{\lceil e/2\rceil}=p^{\lceil \frac{\nu _p(D)}{2}\rceil}$
if $\nu _p(k+1)<\lceil \frac{\nu _p(D)}{2}\rceil$ and $P_{p,k,f}=1$ if
$\nu _p(k+1)\ge \lceil \frac{\nu _p(D)}{2}\rceil$. Thus Lemma 4.6 is true
if either $k\ge p^{\lceil \frac{\nu _p(D)}{2}\rceil}$
and  $2\nmid \nu _p(D)$, or $k\ge p^{\lceil \frac{\nu _p(D)}{2}\rceil}$
and $(\frac{D_p}{p})=-1$.

If $2\mid \nu _p(D)$ and $(\frac{D_p}{p})=1$,
then $\lceil \frac{\nu _p(D)}{2}\rceil=\frac{\nu _p(D)}{2}$ and so
$\nu _p(L_k)\ge \frac{\nu _p(D)}{2}$. First we claim that
\begin{align}\label{4.42}
\nu _p(B_k)=\nu _p(D)/2+\max_{1\le i\le k \atop \nu _p(i)=\nu _p(D)/2}\{ \nu _p(a^2i^2-D)\}.
\end{align}
Let
$
C_1:=\max_{1\le i\le k \atop \nu _p(i)=\nu _p(D)/2}\{ \nu _p(i(a^2i^2-D))\}\
\mbox{and}\ C_2:=\max_{1\le i\le k \atop \nu _p(i)>\nu _p(D)/2}\{ \nu _p(i(a^2i^2-D))\}.
$
Since  $C_1\ge \frac{3\nu _p(D)}{2}$ and $\nu _p(i(a^2i^2-D))<\frac{3\nu _p(D)}{2}$
if $\nu _p(i)<\frac{\nu _p(D)}{2}$, we have by (\ref{4.39}) that
\begin{align}\label{4.43}
\nu _p(B_k)=\max(C_1, C_2).
\end{align}
It also implies that $\nu _p(B_k)\ge \frac{3\nu _p(D)}{2}$.

Note that $C_2=\max_{1\le i\le k \atop \nu _p(i)>\nu _p(D)/2}\{ \nu _p(i)\}+\nu _p(D)=\nu _p(L_k)+\nu _p(D)$
if $\nu _p(L_k)>\frac{\nu _p(D)}{2}$. Thus by (\ref{4.43}), we obtain that
\begin{align}\label{4.44}
\nu _p(B_k)=\max\big( C_1, \nu _p(L_k)+\nu _p(D) \big)
\end{align}
if $\nu _p(L_k)>\frac{\nu _p(D)}{2}$.
If $\nu _p(L_k)>\frac{\nu _p(D)}{2}$, since the discriminant of $a^2x^2-D$ is $4a^2D$
and $\nu _p(4a^2D)=\nu _p(D)$, applying Lemma 3.3 (iii) to the congruence
$a^2x^2-D\equiv0\pmod {p^{\nu _p(L_k)+\frac{\nu _p(D)}{2}}}$, we know that
there is an integer $x_0\in [1, p^{\nu _p(L_k)}]\subseteq [1, k]$ such that
$\nu _p(x_0)=\frac{\nu _p(D)}{2}$ and $\nu _p(a^2x_0^2-D)\ge \nu _p(L_k)+\frac{\nu _p(D)}{2}$.
Hence
$$C_1=\frac{\nu _p(D)}{2}+\max_{1\le i\le k \atop \nu _p(i)=\nu _p(D)/2}
\{ \nu _p(a^2i^2-D)\}\ge \nu _p(L_k)+\nu _p(D).$$
It then follows from (\ref{4.44}) that $\nu _p(B_k)=C_1$ if $\nu _p(L_k)>\frac{\nu _p(D)}{2}$.
On the other hand, there is no integer $i\in [1, k]$ such that
$\nu _p(i)>\frac{\nu _p(D)}{2}$ if $\nu _p(L_k)=\frac{\nu _p(D)}{2}$.
So by (\ref{4.43}), $\nu _p(B_k)=C_1$ if $\nu _p(L_k)=\frac{\nu _p(D)}{2}$.
Thus $\nu _p(B_k)=C_1$ if $k\ge p^{\lceil \frac{\nu _p(D)}{2}\rceil}$
and $2\mid \nu _p(D)$ and $(\frac{D_p}{p})=1$. The claim (\ref{4.42}) is proved.

One can easily check that $\max_{1\le i\le k}\{\nu _p(a^2i^2-D)\}=
\max_{1\le i\le k\atop \nu _p(i)=\nu _p(D)/2}\{\nu _p(a^2i^2-D)\}.$
It then follows from (\ref{4.42}) that
$\nu _p(B_k)=\frac{\nu _p(D)}{2}+\max_{1\le i\le k}\{\nu _p(a^2i^2-D)\}.$
Hence by Lemma 3.13 (ii), we have
$e=\max_{1\le i\le k}\{\nu _p(a^2i^2-D)\}=\nu _p(B_k)-\frac{\nu _p(D)}{2},$
which implies that $e-\frac{\nu _p(D)}{2}=\nu _p(B_k)-\nu _p(D)\ge \frac{\nu _p(D)}{2}$
and so $e\ge \nu _p(D)$. Also we have $\lceil e/2\rceil=\nu _p(D)/2
=e-\nu _p(D)/2=\nu _p(B_k)-\nu _p(D)$ if $e=\nu _p(D)$.
It then follows from Lemma 4.4 that $P_{p,k,f}=p^{\nu _p(B_k)-\nu _p(D)}$
if $\nu _p(k+1)<\nu _p(B_k)-\nu _p(D)$ and $P_{p,k,f}=1$ if $\nu _p(k+1)\ge
\nu _p(B_k)-\nu _p(D)$. So Lemma 4.6 is true
if $k\ge p^{\lceil \frac{\nu _p(D)}{2}\rceil}$, $2\nmid \nu _p(D)$ and
$(\frac{D_p}{p})=1$.
This completes the proof of Lemma 4.6.
\end{proof}

\section{\bf Proof of Theorem \ref{thm1.1} and examples}

In this section, we first give the proof of Theorem \ref{thm1.1}
by using Lemmas 2.2, 4.2, 4.5 and 4.6.

{\it Proof of Theorem \ref{thm1.1}.} Since $g_{k, f}(n)=g_{k, -f}(n)$ for any
$n\in \mathbb{N}^*\setminus Z_{k, f}$, we can assume that $a>0$ in the
following. If $f(x)=ax^2+bx+c$ satisfies that $\gcd(a,b,c)=d>1$,
we can then easily get that
$$g_{k,f}(n)=\frac{\prod_{i=0}^k |f(n+i)|}
{{\rm lcm}_{0\le i\le k}\{f(n+i)\}} =d^k\frac{\prod_{i=0}^k
f_1(n+i)} {{\rm lcm}_{0\le i\le k}\{f_1(n+i)\}}=d^k g_{k, f_1},$$
where $f_1(x)=a_1x^2+b_1x+c_1$ with $a_1=a/d, b_1=b/d$ and
$c_1=c/d$. Obviously, $g_{k,f}$ and $g_{k,f_1}$ have the same
periodicity. If they are both periodic, they have the same
smallest period. That is, we have $P_{k,f}=P_{k,f_1}$.
Therefore we assume that $f(x)$ is a primitive polynomial
(i.e., $\gcd(a,b,c)=1$) in what follows.

Since $f$ is primitive, by Theorem 2.1 we know that the first part
of Theorem \ref{thm1.1} is true. Now we assume that $D\ne a^2i^2$
for all integers $i$ with $1\le i\le k$. Then $g_{k,f}$ is eventually
periodic by Theorem 2.1. Note that $\mathcal {K}_f$ is nonempty and $k\in \mathcal {K}_f$.
In what follows we determine the smallest period $P_{k, f}$ of $g_{k,f}$.
Let $\Delta _{p,k}:=\nu _p(B_k)-\nu _p(P_{p,k,f})$ for any prime $p$.
Since by Lemma 2.2, $P_{p,k,f}\mid p^{\nu _p(B_k)}$ for any prime $p$.
Hence $P_{2,k,f}=1$ if $2\nmid B_k$. So again by Lemma 2.2, we can derive that
\begin{align}\label{5.1}
P_{k,f} =P_{2,k,f}\prod_{p\ne 2, p|B_k}p^{\nu _p(P_{p, k, f})}
=\frac{B_k}{\displaystyle 2^{\Delta _{2,k}}\prod_{p\ne2,
p|B_k} p^{\Delta _{p,k}}}=\frac{B_k}{E_kF_k},
\end{align}
where
\begin{align}\label{5.2}
E_k:=2^{\Delta _{2,k}}\Big(\displaystyle \prod_{p\ne 2, p|\gcd(a, b)}
p^{\Delta _{p,k}}\Big)\Big(\prod_{p\nmid 2a, p|D}p^{\Delta _{p,k}}\Big)
\Big(\prod_{p\nmid 2aD, (\frac{D}{p})=-1}p^{\Delta _{p,k}}\Big)
\end{align}
and
\begin{align}\label{5.3}
F_k:=\Big( \prod_{p|a, p\nmid 2b}p^{\Delta _{p,k}}\Big)
\Big(\prod_{p\nmid 2aD, (\frac{D}{p})=1}p^{\Delta _{p,k}} \Big).
\end{align}

First we treat with $E_k$. If $p=2$, then we get by Lemmas 4.2 and 4.5 that
\begin{align*}
\nu _2(P_{2, k, f})={\left\{
  \begin{array}{rl}
  \nu _2(B_k), & \mbox{if}\ 2|a, 2\nmid b\ \mbox{and}\ \nu _2(k+1)<\nu _2(B_k),\\
  \nu _2(B_k)-2\nu _2(L_k), & \mbox{if}\ 2\nmid a,\  k<2^{\lfloor \frac{\nu _2(D)}
  {2}\rfloor} \ \mbox{and}\ \nu _2(k+1)<\nu _2(L_k),\\
  \lfloor \frac{\nu _2(D)}  {2}\rfloor, & \mbox{if} \  2\nmid a, k\ge
  2^{\lfloor \frac{\nu _2(D)}{2}\rfloor}, D_4\not\equiv1\pmod8 \\
  & \mbox{and}\ \nu _2(k+1)<\lfloor \frac{\nu _2(D)}{2}\rfloor,\\
\nu _2(B_k)-\nu _2(D)-1, & \mbox{if} \ 2\nmid a, k\ge
  2^{\lfloor \frac{\nu _2(D)}{2}\rfloor} \ \mbox{and}\
  D_4\equiv1 \pmod 8,\\
0, & \mbox{otherwise.}
 \end{array}
\right.}
\end{align*}
Thus $2^{\Delta _{2,k}}=\xi_2$ with $\xi _2$ being defined in (\ref{1.3}).

If $p\ne 2$ and $p|\gcd (a, b)$, then by Lemma 4.2, we have $\nu _p(P_{p,k,f})=0$ and so
$\Delta_{p, k}=\nu _p(B_k)$. Hence
\begin{align}\label{5.5}
\prod_{p\ne 2, p|\gcd(a, b)}p^{\Delta _{p,k}}=\prod_{p\ne 2, p|
\gcd(a, b)}p^{\nu _p(B_k)}.
\end{align}

If $p\nmid 2a$ and $p|D$, then using Lemma 4.6, we obtain
$$
\Delta_{p,k}={\left\{
  \begin{array}{rl}
  2\nu _p(L_k), & \mbox{if} \ k<p^{\lceil \frac{\nu _p(D)}
  {2}\rceil} \ \mbox{and}\ \nu _p(k+1)<\nu _p(L_k),\\
  \nu _p(B_k)-\lceil \frac{\nu _p(D)}{2}\rceil, & \mbox{if} \
  k\ge p^{\lceil \frac{\nu _p(D)}{2}\rceil},\
  \nu _p(k+1)<\lceil \frac{\nu _p(D)}{2}\rceil \\ &
  \mbox{and either} \ 2\nmid \nu _p(D) \  \mbox{or}\ (\frac{D_p}{p})=-1,\\
\nu _p(D), &  \mbox{if} \  k\ge
  p^{\lceil \frac{\nu _p(D)}{2}\rceil}, \nu _p(k+1)<\nu _p(B_k)-\nu _p(D),\\
  & 2|\nu _p(D) \ \mbox{and}\ (\frac{D_p}{p})=1,\\
\nu _p(B_k), & \mbox{otherwise.}\\
 \end{array}
\right.}
$$
It follows that
\begin{align}\label{5.6}
\prod_{p\nmid 2a, p|D}p^{\Delta _{p,k}}=\prod_{p\nmid 2a, p|D}\eta_p,
\end{align}
where $\eta_p$ is defined as in (\ref{1.4}).

If $p\nmid 2aD$ and $(\frac{D}{p})=-1$, then $\nu _p(D)=0$ and $D=D_p$.
Hence we have $k\ge p^{\lceil \frac{\nu _p(D)}{2}\rceil}=1$,
$\nu _p(k+1)\ge \lceil \frac{\nu _p(D)}{2}\rceil$ and $(\frac{D_p}{p})=-1$.
It then follows from Lemma 4.6 that $\nu _p(P_{p,k,f})=0$, which implies that
$\Delta_{p, k}=\nu _p(B_k)$. Therefore
\begin{align}\label{5.7}
\prod_{p\nmid 2aD, (\frac{D}{p})=-1}p^{\Delta _{p,k}}
=\prod_{p\nmid 2aD, (\frac{D}{p})=-1}p^{\nu _p(B_k)}.
\end{align}
Now by (\ref{5.2}) and (\ref{5.5})-(\ref{5.7}),
we get that $E_k=\frac{B_k}{A_k}$, where $A_k$
is defined in (\ref{1.2}). Thus by (\ref{5.1}),
we have $P_{k, f}=\frac{A_k}{F_k}.$

Consequently, we handle $F_k$. For this purpose, we first prove
the following fact: There is at most one prime
$p$ such that $\nu _p(k+1)\ge \nu _p(B_k)\ge 1$.
Suppose that there are two distinct primes $p_1$ and $p_2$
such that $v_{p_1}(k+1)\ge v_{p_1}(B_k)\ge 1$
and $v_{p_2}(k+1)\ge v_{p_2}(B_k)\ge 1$. Then $k+1$ is composite and so
$p_1\le k$ and $p_2\le k$. Hence for each $1\le j\le 2$,
$$v_{p_j}(k+1)\ge v_{p_j}(B_k)=\max_{1\le i\le k}\{v_{p_j}(i)+v_{p_j}(a^2i^2-D)\}
\ge \max_{1\le i\le k}\{v_{p_j}(i)\}=v_{p_j}(L_k)\ge 1.$$
But Farhi and Kane \cite{[FK]} showed that there is at most one
prime $p\le k$ such that $\nu _p(k+1)\ge \nu _p(L_k)\ge 1$. We arrive at a
contradiction. Thus the fact is proved.

Now we turn to $F_k$. Let $p\mid B_k$ be a prime satisfying that
either $p|a$ and $p\nmid 2b$ or $p\nmid 2aD$ and $(\frac{D}{p})=1$.
Then $\nu _p(A_k)=\nu _p(B_k)$. It then follows from the above fact that
there is at most one prime $p$ such that $\nu _p(k+1)\ge \nu _p(A_k)\ge 1$.

For any prime $p$ satisfying that either $p|a, p\nmid 2b$
and $\nu _p(k+1)<\nu _p(B_k)$, or $p\nmid 2aD, (\frac{D}{p})=1$
and $\nu _p(k+1)<\nu _p(B_k)$, by Lemmas 4.2 and 4.6,
we deduce that $\nu _p(P_{p, k, f})=\nu _p(B_k)$ and so $\Delta_{p, k}=0$.
If there is no prime $p$ satisfying that $\nu _p(k+1)\ge \nu _p(A_k)\ge 1$
and either $p|a$ and $p\nmid b$ or $p\nmid 2aD$ and $(\frac{D}{p})=1$,
it then follows from (\ref{5.3}) that $F_k=1.$

If there is exactly one odd prime $q$ satisfying that
$v_q(k+1)\ge v_q(A_k)\ge 1$ and either $q|a$ and $q\nmid b$
or $q\nmid 2aD$ and $(\frac{D}{q})=1$, then $v_q(P_{q, k, f})=0$.
Thus $F_k=q^{v_q(B_k)}=q^{v_q(A_k)}.$
Since $P_{k, f}=\frac{A_k}{F_k}$, one concludes that
$P_{k, f}=A_k$ except that $v_q(k+1)\ge v_q(A_k)\ge 1$
for at most one odd prime $q$ such that either $q|a$ and $q\nmid b$
or $q\nmid 2aD$ and $(\frac{D}{q})=1$, in which case one has
$P_{k, f}=A_k/q^{v_q(A_k)}$.
The proof of Theorem \ref{thm1.1} is complete. \hfill$\square$

Now we give two examples to illustrate Theorem \ref{thm1.1}.

\noindent{\bf Example 5.1.} Let $f(x)=4^{l}x^2+1$ with $l\ge 1$
being an integer. Then  $D=-4^{l+1}$,
${\mathcal K}_f=\mathbb{N}^*$ and $Z_{k, f}$ is empty
for all integers $k\ge 1$. By Theorem \ref{thm1.1}, $g_{k, f}$
is periodic for all integers $k\ge 1$.
We have by (1.4), $B_k:={\rm lcm}_{1\le i\le k}\{i(16^{l}i^2+4^{l+1})\}$.
Since $2\mid 4^l=\gcd(a, b)$, we obtain $\xi_2=2^{\nu _2(B_k)}$
by (\ref{1.3}). Clearly, there is no odd prime $p$ such that
$p\mid \gcd(a, b)$
or $p\mid D$. On the other hand, all the primes satisfying
$(\frac{-4^{l+1}}{p})=1$ are of the form $p\equiv1\pmod4$,
and all the primes such that $(\frac{-4^{l+1}}{p})=-1$ must be of the
form $p\equiv3\pmod4$. Hence by (\ref{1.2}), we have
$A_k:=B_k 2^{-\nu _2(B_k)}\prod_{p\equiv3\pmod4}p^{-\nu _p(B_k)}.$
By Theorem \ref{thm1.1}, the smallest period of $g_{k, f}$ equals
$A_k$ except that $\nu _p(k+1)\ge \nu _p(A_k)\ge 1$ for at most
one prime $p\equiv1\pmod{4}$, in which case its smallest
period is equal to $\frac{A_k}{p^{\nu _p(A_k)}}$.

\noindent{\bf Example 5.2.} Let $f(x)=(x+m)(x+m+l)$ with $m\ge 0$
and $l\ge 2$ being integers. Then $D=l^2$ and $\mathcal{K}_f=\{1,..., l-1\}$.
By Theorem \ref{thm1.1}, $g_{k, f}$ is eventually periodic if and only if
$1\le k\le l-1$. Now let $1\le k\le l-1$. Then by (1.4),
$B_k={\rm lcm}_{1\le i\le k}\{ i(a^2i^2-D)\}={\rm lcm}_{1\le i\le k}\{i(l^2-i^2)\}.$
Since $D_4=\frac{l^2}{2^{2\nu _2(l)}}\equiv1\pmod 8$ and $v_q(D)=2v_q(l)$
for any prime factor $q$ of $D$, we have by (\ref{1.3}) and (\ref{1.4}) that

$$\xi_2={\left\{
  \begin{array}{rl}
2^{2\nu _2(L_k)},&   \mbox{if}\ k<2^{\nu _2(l)}\ \mbox{and}\   \nu _2(k+1)<\nu _2(L_k),\\
2^{\nu _2(B_k)}, &  \mbox{if}\ k<2^{\nu _2(l)}\ \mbox{and}\   \nu _2(k+1)\ge \nu _2(L_k),\\
2^{2\nu _2(l)+1}, & \mbox{if}\ k\ge 2^{\nu _2(l)}\\
 \end{array}
\right.}
$$
and
$$
\eta_p={\left\{
  \begin{array}{rl}
p^{2\nu _p(L_k)},&   \mbox{if}\ k<p^{\nu _p(l)}\ \mbox{and}\   \nu _p(k+1)<\nu _p(L_k),\\
p^{2\nu _p(l)}, & \mbox{if}\ k\ge p^{\nu _p(l)}\ \mbox{and}\  \nu _p(k+1)<\nu _p(B_k)-2\nu _p(l),\\
p^{\nu _p(B_k)}, & \mbox{otherwise}.\\
 \end{array}
\right.}$$
Moreover, $(\frac{D}{p})=(\frac{l^2}{p})=1$ for any prime $p$ with $p\nmid 2aD$ and
there is no odd prime $p$ such that $p\mid \gcd(a, b)$. It then follows from (\ref{1.2}) that
$A_k=B_k\xi_2^{-1}\prod_{p\ne 2, p\mid l}\eta_p^{-1}.$
Thus by Theorem \ref{thm1.1}, the smallest period of $g_{k,f}$ is equal to $A_k$ except that
$\nu _p(k+1)\ge \nu _p(A_k)\ge 1$ for at most one odd prime $p$ with $p\nmid 2D$,
in which case its smallest period equals $\frac{A_k}{p^{\nu _p(A_k)}}$.

\section{\bf Proof of Theorem \ref{thm1.2}}

In this section, we show Theorem \ref{thm1.2}.

{\it Proof of Theorem \ref{thm1.2}.} It is clear that if $\gcd(a, b,c)=d$, then
$\log {\rm lcm}_{0\le i\le k}\{f(n+i)\}=
\log {\rm lcm}_{0\le i\le k}\{f_1(n+i)\}+O(1),$
where $f_1(x)= f(x)/d$ is a primitive polynomial. So without loss of generality,
we assume that $\gcd(a, b, c)=1$ and $a>0$ in what follows.

(i). Since $D\ne a^2i^2$ for all $1\le i\le k$,
${\mathcal K}_f$ is nonempty and $k\in {\mathcal K}_f$.
By Theorem \ref{thm1.1}, $g_{k, f}$ is eventually periodic.
So there is a positive integer $n_0$
such that for all positive integers $n\ge n_0$, we have
$g_{k, f}(n)\le M:=\max_{1\le m\le P_{k, f}}\{ g_{k, f}(n_0+m)\}.$
Hence for sufficiently large $n$,
$\log\big(\prod_{i=0}^k|f(n+i)|\big)-\log M\le \log {\rm lcm}_{0\le i\le
k}\{f(n+i)\}\le \log \big(\prod_{i=0}^k|f(n+i)|\big).$

Since $\log \big(\prod_{i=0}^k|f(n+i)|\big)=2(k+1)\log
n+\sum_{i=0}^k\log \Big(a+\frac{2ai+b}{n}+\frac{ai^2+bi+c}{n^2}\Big)$
for sufficiently large $n$, one has
$\lim_{n\rightarrow\infty}\frac{\log
\big(\prod_{i=0}^k|f(n+i)|\big)}{ 2(k+1) \log n}=1 \ {\rm and} \
\lim_{n\rightarrow\infty}\frac{\log
\big(\prod_{i=0}^k|f(n+i)|\big)-\log M}{ 2(k+1) \log n}\\=1.$
This concludes the desired result.

(ii). Since $D=a^2i_0^2$ for some integer $i_0$ with $1\le i_0\le k$,
$f(x)$ is reducible. It then follows from the proof of Theorem 2.1 that $f(x)$ must
be of the form  $(a_1x+b_1)(a_1x+b_1+a_1i_0)$ for some integers
$a_1>0$ and $b_1$ with $\gcd(a_1, b_1)=1$. Thus
${\rm lcm}_{0\le i\le k}\{ f(n+i)\}={\rm lcm}_{0\le i\le k}
\{ \big(a_1(n+i)+b_1\big)\big(a_1(n+i+i_0)+b_1\big)\}.$
It is easy to see that $\{a_1(x+i)+b_1\}_{0\le i\le k+i_0}$
is equal to the set of all the linear factors of $\prod_{i=0}^kf(x+i)$.
Hence ${\rm lcm}_{0\le i\le k+i_0}\{ a_1(n+i)+b_1 \}$ divides
${\rm lcm}_{0\le i\le k}\{ f(n+i)\}$. So we get that
$
{\rm lcm}_{0\le i\le k+i_0}\{ a_1(n+i)+b_1 \}\le {\rm lcm}_{0\le
i\le k}\{ f(n+i)\}\le \prod_{0\le i\le k+i_0}\big(a_1(n+i)+b_1\big)
$
for sufficiently large integer $n$. 

If $b_1\ge 0$, then as in \cite{[HQ]}, we define the following arithmetic function
$g_{k+i_0, a_1, b_1}(n):=\frac{\prod_{0\le i\le
k+i_0}(a_1(n+i)+b_1)}{ {\rm lcm}_{0\le i\le k+i_0}
\{ a_1(n+i)+b_1\}}.$
Then by Theorem \ref{thm1.2} of \cite{[HQ]}, $g_{k+i_0, a_1, b_1}$
is a periodic arithmetic function. So there is a fixed
positive integer $M$ such that $g_{k+i_0, a_1, b_1}(n)\le M$
for all positive integers $n$. If $b_1<0$, then we make a
revision to the above argument by defining
$\tilde g_{k+i_0, a_1, b_1}$ as follows
$\tilde g_{k+i_0, a_1, b_1}(n):=g_{k+i_0, a_1, b_1}(n-b_1).$
Then Theorem \ref{thm1.2} of \cite{[HQ]} tells us that $\tilde g_{k+i_0, a_1, b_1}$
is a periodic arithmetic function. Thus there exists a fixed
positive integer $M$ such that $\tilde g_{k+i_0, a_1, b_1}(n)\le M$
for all positive integers $n$. So $g_{k+i_0, a_1, b_1}(n)\le M$
for all positive integers $n\ge -b_1$.
This concludes that $g_{k+i_0, a_1, b_1}(n)\le M$
for all sufficiently large integers $n$. Thus we obtain that
$\frac{\prod_{0\le i\le k+i_0}(a_1(n+i)+b_1)}{M}\le {\rm lcm}_{0\le
i\le k}\{ f(n+i)\}\le \prod_{0\le i\le k+i_0}(a_1(n+i)+b_1)$
for sufficiently large $n$.
Since $\lim_{n\rightarrow\infty}\frac{\log \prod_{0\le i\le
k+i_0}(a_1(n+i)+b_1)}{(k+i_0+1) \log n}=1,$ the desired result
then follows immediately. So Theorem \ref{thm1.2} is proved. \hfill$\Box$\\
\\
\noindent{\bf Acknowledgements.} The authors would like to thank
Professor Freydoon Shahidi and the anonymous referee for helpful
comments and suggestions which improved its presentation.


\begin{thebibliography}{99}
\bibitem{[A]} T.M. Apostol, Introduction to analytic number
theory, Springer-Verlag, New York, 1976.
\bibitem{[BKS]} P. Bateman, J. Kalb and A. Stenger, A limit involving
least common multiples, Amer. Math. Monthly 109 (2002), 393-394.
\bibitem{[Ch]} P.L. Chebyshev, Memoire sur les nombres premiers, J. Math.
Pures Appl. 17 (1852), 366-390.
\bibitem{[DFI]} W. Duke, J. Friedlander and H. Iwaniec, Equidistribution
of roots of a quadratic congruence to prime moduli, Ann. Math. (2) 141 (1995), 423-441.
\bibitem{[F1]} B. Farhi, Minoration non triviales du plus petit commun multiple de
certaines suites finies d'entiers, C. R. Acad. Sci. Paris S\'er. I 341 (2005), 469-474.
\bibitem{[F2]} B. Farhi, {Nontrivial lower bounds for the least common multiple
of some finite sequences of integers,} J. Number Theory 125 (2007), 393-411.
\bibitem{[FK]} B. Farhi and D. Kane, New results on the least common multiple of
consecutive integers, Proc. Amer. Math. Soc. 137 (2009), 1933-1939.
\bibitem{[Ha]} D. Hanson, {On the product of the primes,} Canad. Math. Bull. 15 (1972), 33-37.
\bibitem{[Ho]} S. Hong, The least common multiple of consecutive terms in a
sequence of positive integers, Proceedings of Sixth International Congress of Chinese
Mathematicians. AMS/IP Stud. Adv. Math., Amer. Math. Soc., Providence, RI, to appear in 2014.
\bibitem{[HF]} S. Hong and W. Feng, {Lower bounds for the least common multiple
of finite arithmetic progressions,} C.R. Acad. Sci. Paris, Ser. I 343 (2006), 695-698.
\bibitem{[HLQW]}S. Hong, Y. Luo, G. Qian and C. Wang, Uniform lower bound for
the least common multiple of a polynomial sequence, C.R. Acad. Sci. Paris,
Ser. I 351 (2013), 781-785.
\bibitem{[HQ]} S. Hong and G. Qian, The least common multiple of consecutive
arithmetic progression terms,  Proc. Edinburgh Math. Soc. 54 (2011), 431-441.
\bibitem{[HQT]} S. Hong, G. Qian and Q. Tan, The least common multiple
of a sequence of products of linear polynomials, Acta Math. Hungar. 135 (2012), 160-167.
\bibitem{[HY]} S. Hong and Y. Yang, On the periodicity of an arithmetical function,
C.R. Acad. Sci. Paris S\'er. I 346 (2008), 717-721.
\bibitem{[HY2]} S. Hong and Y. Yang, Improvements of lower bounds for the least
common multiple of finite arithmetic progressions, {Proc. Amer.
Math. Soc.} 136 (2008), 4111-4114.
\bibitem{[Hu]} L.-K. Hua, {Introduction to number theory,} Springer-Verlag,
Berlin Heidelberg, 1982.
\bibitem{[KK]} D.M. Kane and S.D. Kominers, Asymptotic improvements of
lower bounds for the least common multiples of arithmetic
progressions, Canad. Math. Bull., to appear. http://math.stanford.edu/~dankane/lcmbound.pdf
\bibitem{[K]} N. Koblitz, $p$-Adic numbers, $p$-adic analysis and zeta functions,
GTM 58, Spinger-Verlag, New York, 1984.
\bibitem{[M]} P.J. McCarthy, Introdcution to arithmetical functions,
Springer-Verlag, New York, 1982.
\bibitem{[N]} M. Nair, {On Chebyshev-type inequalities for primes,}
Amer. Math. Monthly 89 (1982), 126-129.
\bibitem{[O]} S.M. Oon, Note on the lower bound of least common multiple,
Abstr. Appl. Anal. 2013, Art. ID 218125, 4 pp.
\bibitem{[QH]} G. Qian and S. Hong, Asymptotic behavior of the least
common multiple of consecutive arithmetic progression terms,
Arch. Math. 100 (2013), 337-345.
\bibitem{[QTH]} G. Qian, Q. Tan and S. Hong, The least common multiple of consecutive
terms in a quadratic progression, {Bull. Aust. Math. Soc.} 86 (2012), 389-404.
\bibitem{[WTH]}R. Wu, Q. Tan and S. Hong, New lower bounds for the least common
multiples of arithmetic progressions, Chin. Ann. Math. 34B (2013), 861-864.
\end{thebibliography}
\end{document}